\documentclass[preprint,final,3p,times]{elsarticle}
\usepackage{lineno}
\usepackage{mathptmx}
\usepackage{amsfonts}
\usepackage{amsmath}
\usepackage{amssymb}
\usepackage[most]{tcolorbox}

\usepackage{hyperref}
\usepackage{xparse}
\usepackage{hhline}
\usepackage[author={Yous van Halder}]{pdfcomment}
\usepackage{lipsum}
\usepackage{wrapfig}
\usepackage{textalpha}

\def\exampletext{Example} 

\NewDocumentEnvironment{testexample}{ O{} }
{
\colorlet{colexam}{black} 
\newtcolorbox[use counter=testexample]{testexamplebox}{%
    empty,
    title={\exampletext: #1},
    attach boxed title to top left,
       minipage boxed title,
    boxed title style={empty,size=minimal,toprule=0pt,top=4pt,left=3mm,overlay={}},
    coltitle=colexam,fonttitle=\bfseries,
    before=\par\medskip\noindent,parbox=false,boxsep=0pt,left=3mm,right=0mm,top=2pt,breakable,pad at break=0mm,
       before upper=\csname @totalleftmargin\endcsname0pt, 
    overlay unbroken={\draw[colexam,line width=.5pt] ([xshift=-0pt]title.north west) -- ([xshift=-0pt]frame.south west); },
    overlay first={\draw[colexam,line width=.5pt] ([xshift=-0pt]title.north west) -- ([xshift=-0pt]frame.south west); },
    overlay middle={\draw[colexam,line width=.5pt] ([xshift=-0pt]frame.north west) -- ([xshift=-0pt]frame.south west); },
    overlay last={\draw[colexam,line width=.5pt] ([xshift=-0pt]frame.north west) -- ([xshift=-0pt]frame.south west); },%
    }
\begin{testexamplebox}}
{\end{testexamplebox}\endlist}

\providecommand{\mbf}[1]{\mathbf{#1}}

\providecommand{\der}[2]{\frac{\partial #1}{\partial #2}}

\newtheorem{theorem}{Theorem}

\newdefinition{rmk}{Remark}
\newproof{pf}{Proof}
\newproof{pot}{Proof of Theorem \ref{thm2}}

\journal{Journal of Computational Physics}

\begin{document}
\begin{frontmatter}

\title{Neural networks for multi-fidelity solver enhancement}

\title{Multi-level neural networks for PDEs with uncertain parameters}

\author[main]{Y. van Halder}
\address[main]{Centrum Wiskunde \& Informatica (CWI), Science Park 123, 1098 XG, Amsterdam, the Netherlands}
\ead{y.van.halder@cwi.nl}
\author[main]{B. Sanderse}
\cortext[cor1]{Corresponding author}
\author[barry]{B. Koren}
\address[barry]{Eindhoven University of Technology, P.O. Box 513, 5600 MB, Eindhoven, the Netherlands}

\begin{abstract}
A novel multi-level method for partial differential equations with uncertain parameters is proposed. The principle behind the method is that the error between grid levels in multi-level methods has a spatial structure that is by good approximation independent of the actual grid level. Our method learns this structure by employing a sequence of convolutional neural networks, that are well-suited to automatically detect local error features as latent quantities of the solution. Furthermore, by using the concept of transfer learning, the information of coarse grid levels is reused on fine grid levels in order to minimize the required number of samples on fine levels. The method outperforms state-of-the-art multi-level methods, especially in the case when complex PDEs (such as single-phase and free-surface flow problems) are concerned, or when high accuracy is required.
\end{abstract}

\end{frontmatter}

\section{Introduction}
\noindent Uncertainty Quantification (UQ) has become increasingly important for complex engineering applications. Determining and quantifying the influence of parametric and  model-form uncertainties is essential for a wide range of applications: from turbulent flow phenomena \cite{xiao_quantifying_2016, edeling_simplex-stochastic_2016}, aerodynamics \cite{simon_gpc-based_2010}, biology \cite{cho_experimentaldesign_2003, abagyan_biased_1994} to design optimisation \cite{constantinescu_computational_2011, mateos_monte_2000, papadrakakis_reliability-based_2002}.

When performing non-intrusive UQ, a black-box is used to sample solutions of a deterministic model, often a PDE, in parameter space, and to interpolate these samples to construct a surrogate model \cite{xiu_numerical_2010}. The number of samples increases exponentially with the number of dimensions of the stochastic space, i.e., the number of uncertainties in the model. This \textit{curse of dimensionality} limits the applicability of UQ algorithms, especially when the black-box is computationally expensive to sample from. Using a low-fidelity black-box alleviates the computational burden by lowering the computational time per sample, but significantly drops the accuracy of the deterministic samples. Using a high-fidelity black-box surely increases accuracy, but at the cost of a significant increase in computational time. As a remedy, multi-fidelity approaches were introduced, which use high-fidelity samples to construct the surrogate model and low-fidelity samples to explore the parameter space, and to see where to place a new high-fidelity sample \cite{peherstorfer_survey_2018, zhu_multi-fidelity_2017,forrester_alexander_i.j_multi-fidelity_2007, haji-ali_multi-index_2016, narayan1, narayan2, teckentrup1, teckentrup2}. These approaches significantly reduce the number of high-fidelity samples required to construct an accurate surrogate model. Multi-level stochastic collocation (MLSC) \cite{teckentrup_multilevel_2015, kouri_multilevel_2014, charrier_finite_2013, harbrecht_multilevel_2013, MLSC_adaptive_1, MLSC_adaptive_2} extends this approach by using a hierarchy of grid resolutions, often referred to as levels, and places samples on each grid level to approximate the difference between two consecutive grid resolutions. The underlying concept of variance decay (of the difference in the solution between subsequent levels) ensures that the number of samples required to approximate the difference between two consecutive levels, reduces with increasing level \cite{MLMC_path}. This results in a significant reduction of required number of high-fidelity samples.

MLSC significantly reduces the required number of high-fidelity samples, due to the utilisation of variance decay of subsequent solution differences. Our method utilises the variance decay and is based on the fact that \textit{subsequent errors themselves typically exhibit similar behaviour (as a function of the spatial variables) with increasing grid resolution}. This is a key insight that we will directly exploit in this paper as it allows us to construct a novel class of multi-level methods based on the concepts of convolutional neural networks and transfer learning. 

Our proposed method is as follows. Like in the MLSC approach \cite{teckentrup_multilevel_2015}, we use the telescopic-sum identity to express the solution in terms of the solution differences between consecutive grid levels ('relative errors', in the following mostly shortly 'errors'). Our method differentiates itself in that we use the similarity of this error between consecutive grid levels to reduce the required number of samples on each level. We achieve this through an original combination of several ingredients. First, we express the relative error between grid levels as a function of the solution at the coarser level through the use of a neural network. Convolutional neural networks enable this step. Convolutional neural networks are well-suited to learn local or low-level features (`latent quantities') in the solution structure, that turn out to be similar to the truncation error commonly encountered in numerical discretisation methods for PDEs. In connection with a second, fully connected, neural network, the convolutional neural network is able to learn the error between grid levels as a function of the solution. Finally, we use transfer learning to reuse the weights and biases from previously trained neural networks at coarse levels to quickly train the neural networks at finer grid levels, thereby significantly reducing the required number of samples of high-fidelity solutions. In the remaining of this paper we refer to the proposed approach as the Multi-Level Neural Network (MLNN) method.


In summary, the highlights of the novel MLNN method for efficient surrogate construction of parametric PDEs are:
\begin{itemize}
\item \textit{The error between grid levels is trained in terms of the solution through a neural network;}
\item \textit{Convolutional neural networks are used to learn local error features as latent quantities;}
\item \textit{Transfer learning is used to exploit the similarities between errors at different levels.}
\end{itemize}

This paper is outlined as follows: section \ref{sec:ProblemDescription} further introduces the problem and notation, section \ref{sec:MF} introduces the machine learning based multi-level approach, section \ref{sec:Method} discusses in detail how to train the neural networks, section 5 summarises the methodology, and finally, section \ref{sec:Results} demonstrates efficiency and accuracy of the MLNN method when applied to several test-cases, ranging from the steady linear advection-diffusion equation to surrogate construction for the unsteady incompressible Navier-Stokes equations.

\section{Problem Description}
\label{sec:ProblemDescription}
\noindent Quantifying the effects of uncertainties in computational engineering typically consists of three steps: (i) the input uncertainties are characterised in terms of a Probability Density Function (PDF), which follows from observations or physical evidence; (ii) the uncertainties are propagated through the model; (iii) the outputs are post-processed, where the Quantity of Interest (QoI) is expressed in terms of its statistical properties. Our goal is to solve the following stochastic problem:
\begin{linenomath*}\begin{equation}
\mathcal{L}(u ; \mbf{z}, \mbf{x}) = 0\ ,\quad \forall \mbf{x}\in D,\quad \forall \mbf{z}\in I\ ,  
\label{eqPD:StochasticProblem}
\end{equation}\end{linenomath*}
where $I\subset \mathbb{R}^k$ is the space of uncertain inputs and is referred to as \textit{random space}, $\mbf{z} = (z_1, ..., z_k)\in I$ the $k$-dimensional vector containing uncertain inputs which follows a joint PDF $\rho(\mbf{z})$, $\mbf{x}\in D$ the vector of spatial/temporal coordinates, $u_{\text{exc}} := u(\mbf{z}, \mbf{x})$ the exact solution, and $\mathcal{L}$ a continuous differential operator. In the remainder of this work we mostly consider $\mbf{x}=(x,y)\in D\subset \mathbb{R}^2$, corresponding to two spatial coordinates. The QoI is given by $q(\mbf{z}, F(u(\mbf{z}, \mbf{x}))$, where $F$ is a given operator.

We are interested in constructing a surrogate model for $q$ as a function of $\mbf{z}$. This is done non-intrusively by sampling solutions $u_i(\mbf{x})$ of \eqref{eqPD:StochasticProblem} at different locations $\{\mbf{z}_i\}_{i=1}^{N_z}$ in the random space and extracting QoI-values from these solutions, i.e.:
\begin{linenomath*}\begin{subequations}\begin{align}
\mathcal{L}(u_i ; \mbf{z}_i, \mbf{x}) &= 0,\quad \forall \mbf{x}\in D,\quad i=1,...,N_z\ ,\label{eqPD:PDE}\\
q_i &= q(\mbf{z_i}) = F(u_i)\ .
\end{align}\end{subequations}\end{linenomath*}
From the set of solutions $\{q_i\}_{i=1}^{N_z}$, a surrogate $\tilde{q}(\mbf{z})$ may be constructed, which satisfies
\begin{linenomath*}\begin{equation}
\tilde{q}(\mbf{z}) \approx q(\mbf{z}),\quad \forall \mbf{z}\in I\ ,
\label{eq:Surrogate}
\end{equation}\end{linenomath*}
and statistical moments, i.e., mean and standard deviation, may be extracted from $\tilde{q}$. However, solving the differential-equation problems \eqref{eqPD:PDE} exactly is often not possible, and therefore the QoI at a given $\mbf{z}_i$ is obtained by solving these problems discretised in space and/or time:
\begin{linenomath*}\begin{subequations}\begin{align}
L(\mbf{u}(\mbf{z}_i); \mbf{z}_i, X) &= 0\ , \label{eqPD:discretisedPDE}\\
q_d(\mbf{z}_i) &= F_d(\mbf{u}(\mbf{z}_i))\ ,
\end{align}\end{subequations}\end{linenomath*}
where the subscript $d$ indicates a quantity which is related to the discretised problem, $L$ the discretised differential operator, $X\in\mathbb{R}^{N\times 2}$ the spatial/temporal computational grid, say $(x_i, y_i)$, consisting of $N$ grid points, $\mbf{u}(\mbf{z}_i)\in\mathbb{R}^{N}$ the solution vector on the computational grid for $\mbf{z}=\mbf{z}_i$, and $F_d$ the operator which maps the discretised solution $\mbf{u}_d(\mbf{z}_i)$ to the QoI $q_d(\mbf{z}_i)$. The computational cost for solving the discretised equation is determined by the discretisation scheme and the resolution of the computational grid $X$, i.e., the number of grid points $N$. If \eqref{eqPD:discretisedPDE} is solved with high precision by using for instance a fine grid or a high-order discretisation scheme, then the solution is referred to as a high-fidelity solution. Contrary, when solving \eqref{eqPD:discretisedPDE} on a coarse grid or by using a low-order discretisation scheme, the solution is referred to as a low-fidelity solution. Other definitions in terms of reduced order models can be given for high/low-fidelity, but they require a different procedure from the one proposed here, and are therefore not considered in this paper. Therefore, the fidelity of the solution is assumed to be fully determined by the grid resolution, and we often refer to different fidelities as levels.

When performing forward propagation of input uncertainties, we have to solve the discretised equation a number of times for different values of $\mbf{z}$. As a result, accurate forward propagation of uncertainties is often infeasible when high-fidelity solutions are computationally expensive to produce. Low-fidelity models may alleviate the computational expenses and make UQ feasible for problems with a complex underlying model, but the accuracy of the solutions decreases and the resulting surrogate model may be inaccurate. Our goal is: \textbf{\textit{to construct a computationally feasible surrogate with high-fidelity accuracy by using solution values of different fidelities}}. The MLNN method uses the main concept of a multi-level method; we use discrete solutions on a hierarchy of grids with increasing resolution. A surrogate is constructed, which is mainly based on a set of low-fidelity solutions, and uses a relatively low number of higher-fidelity solutions in order to enhance the accuracy of the surrogate. The MLNN method is discussed in more detail in the next two sections.

\section{Machine Learning Based Multi-Level Approach}
\label{sec:MF}
\noindent The multi-level approach is explained in this section. To begin with, we introduce notation for the different levels of solutions. The discretised operators corresponding to different levels are labelled with a superscript and are denoted as $L^{(i)}$, and the corresponding computational grids $X^{(i)}$ consisting of $N^{(i)}$ grid points. The solutions are denoted as $\mbf{u}^{(i)}(\mbf{z})\in \mathbb{R}^{N^{(i)}}$ and are functions of the input uncertainties. In the MLNN method we use a hierarchy of grid resolutions consisting of a total of $N_L$ grids. We denote with $\mbf{u}^{(i)}(\mbf{z})$ the solution which is computed on the coarsest grid level $X^{(1)}$, while $\mbf{u}^{(N_L)}(\mbf{z})$ is the solution computed on the finest grid level $X^{(N_L)}$. For simplicity the computational grids are assumed to be nested:
\begin{linenomath*}\begin{equation}
X^{(1)}\subset X^{(2)}\subset ... \subset X^{(N_L)}\ .
\label{eq:nestedness}
\end{equation}\end{linenomath*}

Several multi-level strategies have been presented that use a hierarchy in grid resolutions \cite{peherstorfer_survey_2018}, which reduces the number of high-fidelity solves required to construct an accurate surrogate model or extract accurate statistics. Consider restriction of the fine-grid solution to the coarsest grid:
\begin{linenomath*}\begin{equation}
\mbf{u}^{(N_L)}(\mbf{z})\rightarrow\mbf{u}^{(N_L)}(\mbf{z})|_{X^{(1)}}\ ,
\end{equation}\end{linenomath*}
where $\cdot |_{X^{(1)}}$ is the restriction operator which computes values on the coarse grid $X^{(1)}$, which may incorporate sub sampling and interpolation. The multi-level methods described in \cite{teckentrup_multilevel_2015,charrier_finite_2013,harbrecht_multilevel_2013} rewrite $\mbf{u}^{(N_L)}(\mbf{z})$, using the telescoping-sum identity:
\begin{linenomath*}\begin{equation}
\mbf{u}^{(N_L)}(\mbf{z})|_{X^{(1)}} = \mbf{u}^{(1)}(\mbf{z}) + \sum_{i=2}^{N_L}(\mbf{u}^{(i)}(\mbf{z})|_{X^{(1)}} - \mbf{u}^{(i-1)}(\mbf{z})|_{X^{(1)}})\ ,
\label{eq:telescopic}
\end{equation}\end{linenomath*}
In case of nested grids, \eqref{eq:nestedness}, the restriction operator solely uses sub sampling. Restricting to the coarsest level is not necessary, as we can also prolongate the coarse solutions to the finest level computational grid $X^{(N_L)}$, but in many cases the coarsest grid suffices to accurately calculate the QoI. To circumvent interpolation of coarse-grid solution values to the fine grid and reduce the degrees of freedom in the neural networks that are used later, we opt to restrict solution values to the coarsest level, but the MLNN method can be extended to account for progolongation to the finest level as well. For convenience, we denote the error between two levels as
\begin{linenomath*}\begin{equation}
\mbf{e}^{(i)}(\mbf{z}):= \mbf{u}^{(i)}(\mbf{z})|_{X^{(1)}} - \mbf{u}^{(i-1)}(\mbf{z})|_{X^{(1)}}\ .
\end{equation}\end{linenomath*}
The goal of a multi-level method is to approximate the error between two levels $\mbf{e}^{(i)}$, where we assume that $\mbf{e}^{(i)}\rightarrow 0$ as $i\rightarrow\infty$. As a result, the variance decay in the magnitude of the terms in \eqref{eq:telescopic} decreases the number of samples required to approximate $\mbf{e}^{(i)}(\mbf{z})$ with increasing level $i$ \cite{teckentrup_multilevel_2015, giles2015multilevel}, which reduces the total computational cost for constructing an accurate surrogate. 

The MLNN method is also based on identity \eqref{eq:telescopic}, but does not use interpolation for constructing the approximations for $\mbf{e}^{(i)}$ as used in \cite{teckentrup_multilevel_2015}. Instead we use a convolutional neural network that approximates $\mbf{e}^{(i)}$ in terms of a set of local low-level features, also known as latent quantities. To show that a representation of $\mbf{e}^{(i)}$ in a set of latent quantities is justified, the error $\mbf{e}^{(i)}$ can be written as:
\begin{linenomath*}\begin{align}
\mbf{e}^{(i)} &= \mbf{u}^{(i)}(\mbf{z})|_{X^{(1)}} - \mbf{u}^{(i-1)}(\mbf{z})|_{X^{(1)}}\nonumber\\
&= (\mbf{u}^{(i)}(\mbf{z})|_{X^{(1)}} - \mbf{u}_{\text{exc}}(\mbf{z})|_{X^{(1)}}) - (\mbf{u}^{(i-1)}(\mbf{z})|_{X^{(1)}} - \mbf{u}_{\text{exc}}(\mbf{z})|_{X^{(1)}})\nonumber\\
&= \mbf{E}^{(i)}|_{X^{(1)}} - \mbf{E}^{(i-1)}|_{X^{(1)}}\ ,\label{eq:globalerror}
\end{align}\end{linenomath*}
where $\mbf{E}^{(i)}=\mbf{u}^{(i)}(\mbf{z}) - \mbf{u}_{\text{exc}}(\mbf{z})|_{X^{(i)}}$ denotes the global error between the solution on level $i$ and the exact solution (sub sampled on $X^{(i)}$). The exact solution is often not known and therefore \eqref{eq:globalerror} is not directly useful. However, for a linear discretisation operator $L^{(i)}$, we can write:
\begin{linenomath*}\begin{align}
L^{(i)} \mbf{E}^{(i)} &= L^{(i)} (\mbf{u}^{(i)} - \mbf{u}_{\text{exc}}|_{X^{(i)}}) \nonumber \\
&= L^{(i)}\mbf{u}^{(i)} - L^{(i)} (\mbf{u}_{\text{exc}}|_{X^{(i)}}) \nonumber\\
&=  - L^{(i)} (\mbf{u}_{\text{exc}}|_{X^{(i)}})\nonumber \\
&= -\mbf{\mbf{\tau}}^{(i)}\ ,\label{eq:truncation0}
\end{align}\end{linenomath*}
where $\tau$ is the local truncation error. This can be rewritten as:
\begin{linenomath*}\begin{align}
\mbf{E}^{(i)} &= -\left(L^{(i)}\right)^{-1}\mbf{\tau}^{(i)}\ ,\label{eq:globaltrucation}
\end{align}\end{linenomath*}
The local truncation error is the error when substituting the exact solution into the discretised operator $L^{(i)}$ and is computed locally, for instance by performing Taylor-series expansions on the schemes that are used to discretise the differential operator $\mathcal{L}$. The inverse operator $\left(L^{(i)}\right)^{-1}$ maps the local truncation error to the global error $E^{(i)}$. For linear discretisation of linear PDEs this results in the form:
\begin{linenomath*}\begin{equation}
\mbf{\tau}^{(i)} = \sum_{k=1}^{\infty}c_k \left(\mbf{u}_{\text{exc}}\right)_{\partial, k}(\mbf{z})|_{X^{(i)}}\ ,
\label{eq:truncation}
\end{equation}\end{linenomath*}
where $\left(\mbf{u}_{\text{exc}}\right)_{\partial, k}$ is a sequence of partial derivatives of the exact solution with respect to spatial/temporal coordinates. The truncation error can be interpreted as a set of latent quantities, which are local or low-level features that can be used to effectively express the error in the solution. However, this only holds for linear differential equations and is therefore a rather restrictive assumption.

Our key idea is to construct an approximation for $\mbf{e}^{(i)}$ in terms of a set of latent quantities. This set of latent quantities is obtained by applying a non-linear transformation of the solution values on a coarser level. For linear differential equations, the truncation error is a good candidate for such a set of latent quantities, but it does not generalise to non-linear problems. Therefore, we propose a new algorithm that approximates $\mbf{e}^{(i)}(\mbf{z}),\ i=2,...,N_L$ using a neural network approach. The \textit{How} and  \textit{Why} of this approach are discussed in more detail next.\\

\begin{figure}[hbt!]
\includegraphics[width = \textwidth]{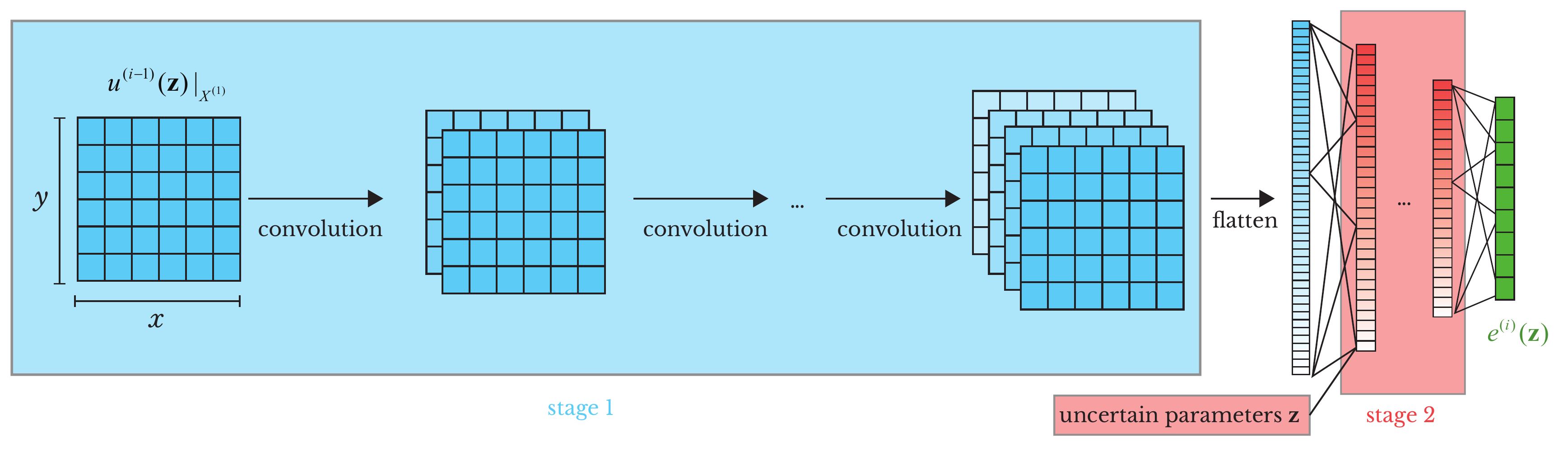}
\caption{\label{fig:MethodOverview} The convolutional and fully-connected network architecture for approximating $\mbf{e}^{(i)}$ for a scalar partial differential equation with two space dimensions and uncertainties $\mbf{z}$.}
\end{figure}

\noindent\textit{How a neural network is used to approximate $\mbf{e}^{(i)}(\mbf{z})$}\\
The neural network architecture that is used for approximating $\mbf{e}^{(i)}(\mbf{z})$ is shown in figure \ref{fig:MethodOverview}. The neural network consists of two stages:\\ \

\indent\textit{Stage 1: Obtain latent/inferred quantities from the solution vector $\mbf{u}^{(i)}(\mbf{z})$, that can be mapped efficiently to $\mbf{e}^{(i)}(\mbf{z})$}\\
\indent\textit{Stage 2: Construct mapping from latent quantities to $\mbf{e}^{(i)}(\mbf{z})$}\\ 

\noindent The neural network takes as input the sub-sampled solution values $\mbf{u}^{(i-1)}(\mbf{z})|_{X^{(1)}}$ and is trained such that it outputs $\mbf{e}^{(i)}(\mbf{z})$. \\

\noindent\textit{Stage 1, How?}\\
The first part of the neural network that corresponds to stage 1 consists of a number of sequential convolutional layers. The filters that are used in convolutional neural networks are equivalent to a sequence of local finite difference operators whose coefficients are determined during the training procedure. As a result, the convolutional layers are used to efficiently combine the solution values into latent quantities.\\

\noindent\textit{Stage 1, Why?}\\
Instead of constructing an approximation for $\mbf{e}^{(i)}(\mbf{z})$ directly (as done for example in \cite{teckentrup_multilevel_2015}), we first construct a set of latent quantities from the solution vector $\mbf{u}^{(i)}(\mbf{z})|_{X^{(1)}})$, which are then mapped to $\mbf{e}^{(i)}(\mbf{z})$. We learn from equations \eqref{eq:truncation0}-\eqref{eq:globaltrucation}, that a mapping between the solution vector and $\mbf{e}^{(i)}(\mbf{z})$ exists in terms of the truncation error, which consists of a set of latent quantities. As a result, the MLNN method uses a sequence of convolutional layers to find a set of latent quantities that can be used by stage 2 of the neural network to efficiently construct an approximation for $\mbf{e}^{(i)}(\mbf{z})$.\\

\noindent\textit{Stage 2, How?}\\
After propagating the input $\mbf{u}^{(i-1)}(\mbf{z})|_{X^{(1)}}$ through the convolutional layers, the output is flattened into a single vector and feeds the second stage of the neural network. The second stage consists of multiple fully-connected layers that apply a non-linear transformation that maps the previously obtained latent quantities to the desired quantity $\mbf{e}^{(i)}(\mbf{z})$. Furthermore, as we want to approximate the dependency of $\mbf{e}^{(i)}$ on $\mbf{z}$, $\mbf{z}$ is included in the neural network by concatenating the values of $\mbf{z}$ to the flattened convolutional output.\\

\noindent\textit{Stage 2, Why?}\\
Solely using convolutional layers is not preferred due to their inherent local nature, making them only useful for extracting low-level features. The latent quantities obtained in stage 1 need to be mapped to the desired quantity  $\mbf{e}^{(i)}(\mbf{z})$. As in equation (12), this requires a global and possibly non-linear transformation. It is known that fully-connected neural networks are so-called universal-function approximators \cite{HORNIK1991251, unversal-nn, hanin2017approximating} and are therefore perfectly suited for constructing the non-linear mapping between the latent quantities and $\mbf{e}^{(i)}(\mbf{z})$.\\

\noindent\textit{Stage 1 + Stage 2}\\
The full architecture for a scalar partial differential equation with two space dimensions and uncertainties $\mbf{z}$ is shown schematically in figure \ref{fig:MethodOverview}.
\\
The architecture for the coupled vector partial differential equation case requires a multi-channel input, i.e., one channel for each of the calculated quantities, which is explained in more detail in section 4. We denote the neural network that maps $\mbf{u}^{(i-1)}(\mbf{z})|_{X^{(1)}}$ to $\mbf{e}^{(i)}(\mbf{z})$ as,
\begin{linenomath*}\begin{equation}
\mbf{e}^{(i)}(\mbf{z}) \approx P^{(i)}(\mbf{u}^{(i-1)}(\mbf{z})|_{X^{(1)}})\ .
\end{equation}\end{linenomath*}
To clarify, the neural network $P^{(i)}$ requires a training set that comprises pairs of solutions $(\mbf{u}^{(i-1)}(\mbf{z})|_{X^{(1)}}, \mbf{u}^{(i)}(\mbf{z})|_{X^{(1)}})$, where $\mbf{u}^{(i-1)}(\mbf{z})|_{X^{(1)}}$ is used as input for the neural network, while $\mbf{u}^{(i)}(\mbf{z})|_{X^{(1)}}$ enters the cost function during training. This is discussed in more detail in section \ref{sec:Method}.
In order to utilise the telescopic-sum identity \eqref{eq:telescopic}, we need approximations for $\mbf{e}^{(i)}(\mbf{z})$ with $i=2,...,N_L$ and therefore we train $N_L-1$ neural networks $(P^{(2)},...,P^{(N_L)})$ of the form shown in figure \ref{fig:MethodOverview}. After training the neural networks, we can approximate the high-fidelity solution as:
\begin{linenomath*}\begin{equation}
\mbf{u}^{(N_L)}(\mbf{z})|_{X^{(1)}} \approx \mbf{u}^{(1)}(\mbf{z}) + \sum_{i=2}^{N_L}P^{(i)}\left(\mbf{u}^{(1)}(\mbf{z}) + \sum_{j=2}^{i-1}P^{(j)}\left(\mbf{u}^{(1)}(\mbf{z}) + \sum_{k=2}^{j-1}P^{(k)}\left( ... \right)\right)\right)\ ,
\label{eq:approximatesolution}
\end{equation}\end{linenomath*}
where we used \eqref{eq:telescopic} and recursively used: 
\begin{linenomath*}\begin{equation}
\mbf{u}^{(i)}(\mbf{z})|_{X^{(1)}} = \mbf{u}^{(i-1)}(\mbf{z})|_{X^{(1)}} + \mbf{e}^{(i)}\approx \mbf{u}^{(i-1)}(\mbf{z})|_{X^{(1)}} + P^{(i)}(\mbf{u}^{(i-1)}(\mbf{z})|_{X^{(1)}})\ .
\end{equation}\end{linenomath*}
We denote the approximate high-fidelity solution as $\tilde{\mbf{u}}^{(N_L)}(\mbf{z})\approx\mbf{u}^{(N_L)}(\mbf{z})|_{X^{(1)}}$. Notice that we only need $\mbf{u}^{(1)}(\mbf{z})|_{X^{(1)}}$ to construct an approximation for $\mbf{u}^{(N_L)}(\mbf{z})|_{X^{(1)}}$. The architecture, hyperparameters and training procedure of the neural networks will be discussed in more detail in section \ref{sec:Method}.\\

\noindent\textit{How the required number of samples for training $P^{(i)}$ is reduced}\\
The combination of convolutional layers and fully-connected layers allows for learning complex non-linear mappings. The filter coefficients that are determined during training of the convolutional layers, have significantly less degrees of freedom when compared to fully-connected layers with similar input-output dimensions. As a result, the number of unknowns that need to be determined during training is significantly lower than in a network architecture with the same number of layers using solely fully-connected layers. This drastically decreases the required number of training samples in the MLNN method.

Second, it can be shown that $\mbf{e}^{(i-1)}(\mbf{z})$ and $\mbf{e}^{(i)}(\mbf{z})$ have a similar spatial shape (apart from a scaling factor). This concept is summarised in the following theorem:
\begin{theorem}
Assume the solution $\mbf{u}$ of the discretised problem \eqref{eqPD:discretisedPDE} can be expanded as follows:
\begin{linenomath*}\begin{equation}
\mbf{u}^{(h)} = \mbf{u}_{\text{exc}}|_{X^{h}} + h^d \mbf{v}|_{X^{h}} + O(h^{d+1})\ ,
\label{eq:expansion}
\end{equation}\end{linenomath*}
where $h$ is the grid resolution, $d$ is the order of convergence with $d\geq 1$, and $v$ is a function independent of $h$. Then the errors $\mbf{e}^{(h/2)}:= \mbf{u}^{(h/2)}|_{X^h} - \mbf{u}^{(h)}$ and $\mbf{e}^{(h/4)}:= \mbf{u}^{(h/4)}|_{X^h} - \mbf{u}^{(h/2)}|_{X^h}$ satisfy
\begin{linenomath*}\begin{equation}
c\mbf{e}^{(h/4)} = \mbf{e}^{(h/2)} + O(h^k)\ , 
\label{eq:wantedsimilar}
\end{equation}\end{linenomath*}
where $k\geq 1$ and $c$ is a constant independent of $h$.
\end{theorem}
\begin{pf}
Using the expansion \eqref{eq:expansion}, we can write
\begin{linenomath*}\begin{align}
\mbf{e}^{(h/2)} = \left(\left(\frac{h}{2}\right)^d - h^d\right)\mbf{v}|_{X^h} + O(h^{d+1})\ ,\\
\mbf{e}^{(h/4)} = \left(\left(\frac{h}{4}\right)^d - \left(\frac{h}{2}\right)^d\right)\mbf{v}|_{X^h} + O(h^{d+1})\ .
\end{align}\end{linenomath*}
Multiplying the second equation with $2^d$ results in
\begin{linenomath*}\begin{equation}
2^d\mbf{e}^{(h/4)} = \left(\left(\frac{h}{2}\right)^d - h^d\right)\mbf{v}|_{X^h} + 2^dO(h^{d+1})  = \mbf{e}^{(h/2)} + O(h^{d+1})\ ,\\
\end{equation}\end{linenomath*}
which is of the form displayed in equation \eqref{eq:wantedsimilar} with $k=d+1$.\hspace*{0pt}\hfill $\square$
\end{pf}
The existence of expansion \eqref{eq:expansion} is a crucial assumption in this theorem. It has been proven in \cite{proofbook} (theorem 2.1) that this expansion exists for a large class of (discretisations of) linear differential equations under relatively mild conditions. Roughly speaking, these conditions are: existence and uniqueness of the discrete and continuous problem, and sufficient smoothness of the terms appearing in the discretised equations. 

Theorem 1 indicates that if an expansion of the form \eqref{eq:expansion} holds, then information of a previously trained neural network mapping $P^{(i-1)}$ can be used to more efficiently train $P^{(i)}$. The concept of reusing parts of a previously trained neural network is not new in the field of machine learning, and is known as transfer learning \cite{aggarwal_neural_2018}. Transfer learning is used to circumvent the need for sampling many solutions on the higher levels by using weights and biases from the previously trained neural network $P^{(i-1)}$. 

The concept of similarity of error profiles, as expressed by theorem 1, only holds for linear differential equations, as it relies on assumption \eqref{eq:expansion}. However, we numerically show in section \ref{sec:Results} that the concept of similarity of error profiles also holds to a large extent for non-linear differential equations, which makes transfer learning applicable for non-linear cases as well.

\section{Training the Neural Networks}
\label{sec:Method}
\noindent In order to train the neural networks $(P^{(2)},...,P^{(N_L)})$, we need to specify:
\begin{enumerate}[I]
\item \textbf{Architecture}: layer architecture, activation functions, and hyperparameters
\item \textbf{Training/validation data}: data set for training and validating the neural networks
\item \textbf{Hyperparameter tuning}: finding proper hyperparameter values
\end{enumerate}
We will discuss these three components separately in the next subsections.

\subsection{Neural Network Architecture}
\label{sec:architecture}
\noindent As mentioned before we use a combination of convolutional layers and fully-connected layers, and the architecture is split in two stages.

\subsubsection{Convolutional Layers}
\begin{figure}[hbt!]
\includegraphics[width = \textwidth]{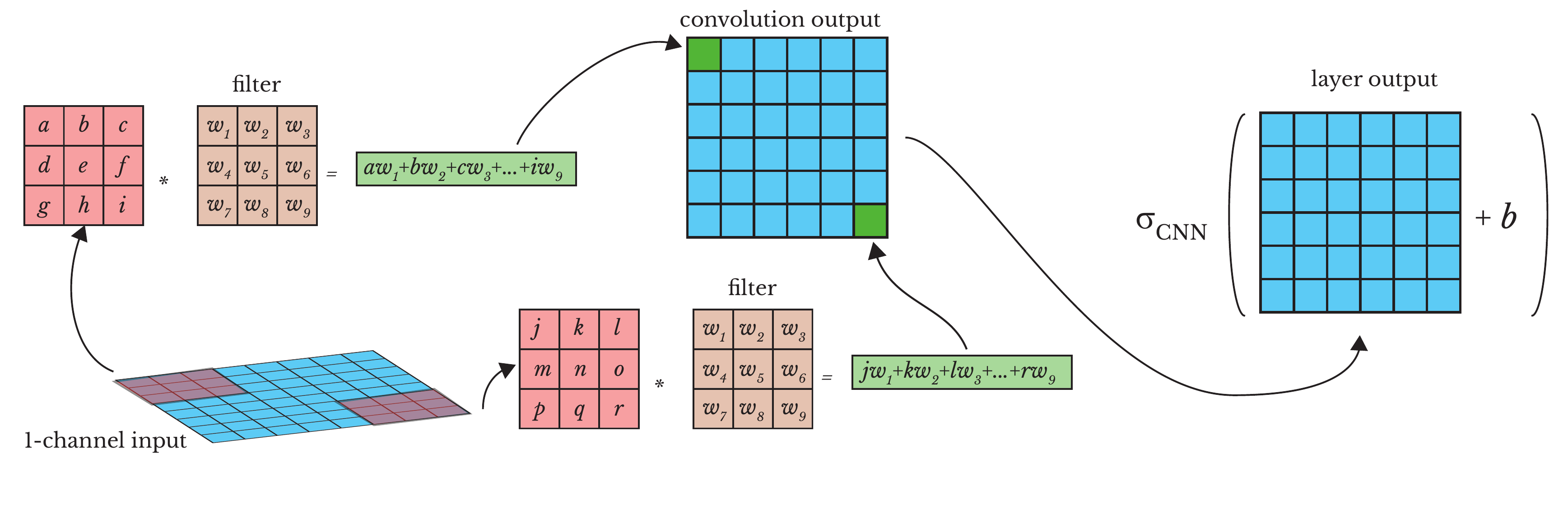}
\caption{\label{fig:CNNlayer} The 2D convolutional layer with a 1-channel input.}
\end{figure}

\noindent The name ``convolutional layer'' indicates that the network employs a mathematical operation called convolution. A convolutional layer is simply a fully-connected neural network layer in which the dense matrix multiplication is replaced by a convolution operator (multiplication with a sparse matrix) \cite{goodfellow_deep_2016}, where filters are convolved over the space of current values. In general, an $n$D convolutional layer accepts $n_c$ channels of dimension $k_1\times k_2 \times ... \times k_n$ as input, and convolves a filter of dimension $f_1\times f_2 \times ... \times f_n \times n_c$ over the input channels, adds a bias, and applies a (non-)linear transformation to produce the output. This procedure is shown schematically in figure \ref{fig:CNNlayer}.
\\
Figure \ref{fig:CNNlayer} shows a 2D convolutional layer for a 1-channel input. The number of channels of the input corresponds to the number of quantities that is solved for in the underlying PDE. E.g., when solving 2D incompressible Navier-Stokes, we solve for the two velocity components and the pressure, which results in a 3-channel input. In this work, a single 3x3 filter is used with 9 unknown filter weights $w_j$, which is a commonly used filter-width that yields the best results on our test cases (see section \ref{sec:Results}). In addition, a multi-channel input results in a multi-channel filter with more unknowns. It is common practice to use more than one filter, which results in an output with multiple channels, i.e., one channel associated to the convolution of one filter. Notice that the output of a convolutional layer is in general smaller than the input, as the filter cannot cross the boundaries of the input channel. Therefore, the dimensions of the output depend on the input dimensions and the filter dimensions. As a remedy, we can pad the input channel with zero values at the boundaries in order to obtain an output which has the same dimensions as the input channel. This is often referred to as zero-padding, and will be used in this paper to keep the output dimensions of the convolutional layers the same. After convolving the filters, a trainable bias coefficient is added to each output channel and the resulting values are activated with a non-linear activation function $\sigma_{\text{CNN}}$. The activation function $\sigma$ for the convolutional neural network (CNN) is chosen to be the Rectified Linear Unit (ReLU):
\begin{linenomath*}\begin{equation}
\sigma_{\text{CNN}}(x) = \max(0, x)\ ,
\end{equation}\end{linenomath*}
because it does not suffer from the vanishing gradient problem and allows for a sparse representation of the output \cite{goodfellow_deep_2016}. The ReLU activation function may suffer from so-called dying neurons \cite{goodfellow_deep_2016}. Other variants, e.g., leaky-Relu or ELU, may be used to resolve the dying neuron problem. However, the sparse representation when using ReLU is beneficial and is therefore employed in the remainder of this work.
The number of convolutional layers $N_{\text{CNN}}$ to choose, depends on the intrinsic complexity of the dataset and is treated as a hyperparameter in the MLNN method. Lastly, increasing the number of filters with a factor 2 for consecutive CNN layers is a good starting point to tune the neural network \cite{goodfellow_deep_2016}, and this strategy is therefore used in this paper.

\subsubsection{Fully-Connected Layers}
\noindent After the input has been propagated through the convolutional layers, we flatten the output of the last convolutional layer into a single vector, and the uncertainties $\mbf{z}$ are concatenated. This vector is the start of the fully-connected part, where we apply a non-linear transformation to the output of the convolutional part. This part of the neural architecture is characterised by the number of fully-connected layers $N_{\text{FC}}$, the number of neurons per layer $n$, and the activation function $\sigma_{\text{FC}}$. The parameters $N_{\text{FC}}$ and $n$ are treated as hyperparameters and will be tuned during the training process. Additionally, a ReLU activation function is used for the fully-connected layers. Once we have propagated the output of the convolutional part through the $N_{\text{FC}}$ fully-connected layers, we output a vector which has the same dimensions as the input of the network. The output layer uses a linear activation function, because it allows for unbounded output values.

\subsubsection{Cost Function}
\noindent The goal of training the neural network is to find the weights $\mbf{w}$ and biases $\mbf{b}$, of both the convolutional and fully-connected parts, which minimise the following cost functions:
\begin{linenomath*}\begin{equation}
c^{(i)}(\mbf{w}, \mbf{b}) = \sum_{(\mbf{u}^{(i)}, \mbf{u}^{(i-1)})\in T^{(i)}} \|P^{(i)}(\mbf{u}^{(i-1)}|_{X^{(1)}}) - (\mbf{u}^{(i)}|_{X^{(1)}} - \mbf{u}^{(i-1)}|_{X^{(1)}})\|_2^2\ ,\ \ i=2,...,N_L\ ,
\label{eq:costfunction}
\end{equation}\end{linenomath*}
where $T^{(i)}$ is the training set. The construction of this training set is discussed in section \ref{sec:trainingdata} in more detail.

Overfitting may occur due to the large amount of unknowns, i.e., convolutional weights/biases and fully-connected layer weights/biases. In order to prevent overfitting, we use an $l_2$-regularisation, which introduces a new hyperparameter $\lambda$ and the altered cost functions:
\begin{linenomath*}\begin{equation}
c_{\lambda}^{(i)}(\mbf{w}, \mbf{b}) = \sum_{(\mbf{u}^{(i)}, \mbf{u}^{(i-1)})\in T} \|P^{(i)}(\mbf{u}^{(i-1)}|_{X^{(1)}}) - (\mbf{u}^{(i)}|_{X^{(1)}} - \mbf{u}^{(i-1)}|_{X^{(1)}})\|_2^2 + \lambda \|\mbf{w}\|_2^2,\quad i=2,...,N_L\ .
\label{eq:costfunctionlambda}
\end{equation}\end{linenomath*}
After training the neural network, we validate if the neural network generalises to unseen data by computing the validation error:
\begin{linenomath*}\begin{equation}
v^{(i)} = \sum_{(\mbf{u}^{(i)}, \mbf{u}^{(i-1)})\in V^{(i)}} \|P^{(i)}(\mbf{u}^{(i-1)}) - (\mbf{u}^{(i)} - \mbf{u}^{(i-1)})\|_2^2,\quad i=2,...,N_L\ ,
\label{eq:validationerror}
\end{equation}\end{linenomath*}
where $V^{(i)}$ is the validation set. The validation error is used to tune the hyperparameters of the network architecture. Additionally, a test set can be used to check if the neural network with proper hyperparameters generalises well outside the training/validation set. However, using a test set results in an increase of required number of samples, and as generalisation is only important within or close to the range of training/validation data, we choose to solely use a training and validation set.

\subsubsection{Transfer Learning}
\begin{figure}[hbt!]
\includegraphics[width = \textwidth]{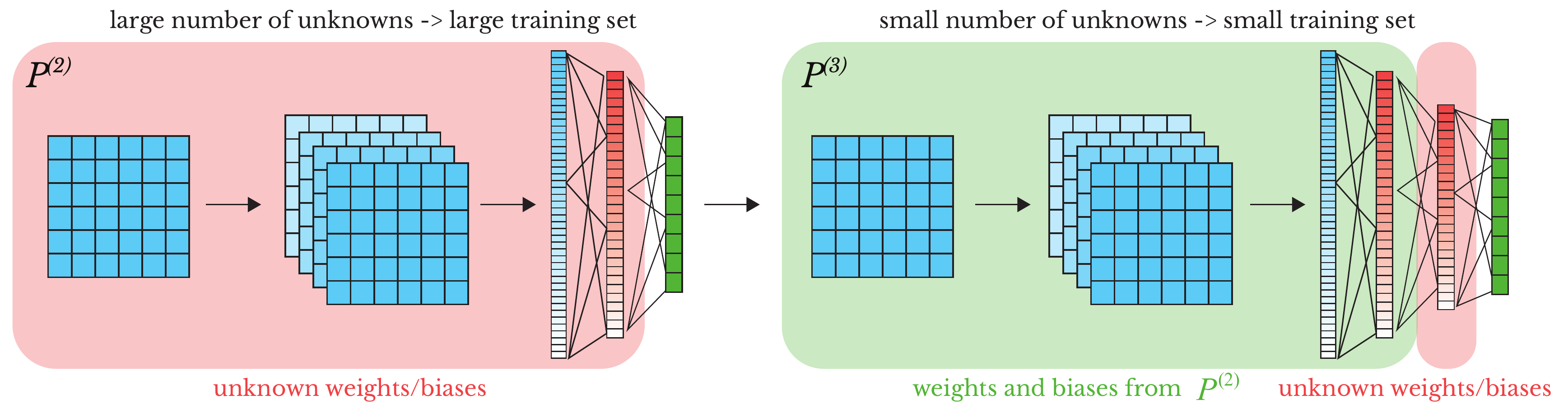}
\caption{\label{fig:transferlearning} Transfer learning procedure.}
\end{figure}

\noindent The amount of training data required for obtaining a proper set of weights and biases may become infeasible when increasing the level, as it requires a training set that is constructed by sampling many high-fidelity solutions. When training $P^{(i)},\ i=3,...,N_L$, transfer learning is used to circumvent the need for sampling many solutions on the higher levels, by using weights and biases from the previously trained neural network $P^{(i-1)}$. Transfer learning can be performed in multiple ways \cite{aggarwal_neural_2018}. We choose to fix the weights and biases of the previously trained neural network $P^{(i-1)}$ and add a small number of fully-connected layers at the end, while keeping the rest of the architecture the same. As the input for $P^{(i)}, i=2,...,N_L$ is a solution vector, which is defined on the same grid $X^{(1)}$, the learned convolutional filters are still expected to compute a proper set of latent quantities when increasing the level $i$. Therefore, the transfer learning approach leaves the convolutional part unaltered for all $P^{(i)},\ i=3,...,N_L$. The transfer learning procedure is schematically shown in figure \ref{fig:transferlearning}.
A single fully-connected layer is added with the same activation function as the other fully-connected layers, but with a hyperparameter $n$ representing the number of neurons in the new layer.

\subsection{Training/Validation Data}
\label{sec:trainingdata}
\noindent The neural networks $\{P^{(i)}\}_{i=2}^{N_L}$ are trained to approximate the errors $\mbf{e}^{(i)}(\mbf{z})|_{X^{(1)}} = \mbf{u}^{(i)}(\mbf{z})|_{X^{(1)}} - \mbf{u}^{(i-1)}(\mbf{z})|_{X^{(1)}}$. The training set for training $P^{(i)}$ requires solution values for both $\mbf{u}^{(i-1)}(\mbf{z})$ and $\mbf{u}^{(i)}(\mbf{z})$, sampled at multiple values for $\mbf{z}$. The training set that is used for training $P^{(i)}$ is denoted as:
\begin{linenomath*}\begin{equation}
T^{(i)}:=\left\lbrace \left(\mbf{u}^{(i)}(\mbf{z_k})|_{X^{(1)}}, \mbf{u}^{(i-1)}(\mbf{z_k})|_{X^{(1)}}\right)\ |\ k=1,...,N_{\text{training}}\right\rbrace\ .
\label{eq:trainingset}
\end{equation}\end{linenomath*}
In order to validate if the neural network generalises to values not in the training set, we validate our neural network on the set:
\begin{linenomath*}\begin{equation}
V^{(i)}:=\left\lbrace \left(\mbf{u}^{(i)}(\mbf{z_k})|_{X^{(1)}}, \mbf{u}^{(i-1)}(\mbf{z_k})|_{X^{(1)}}\right)\ |\ k=1,...,N_{\text{validate}}\right\rbrace\ .
\label{eq:trainingset}
\end{equation}\end{linenomath*}

The amount and quality of data that is used to train and validate the neural network is paramount. How to ensure that we obtain a proper set of weights and biases that generalises well to unseen values, is discussed in this section. The required amount of training data depends on the complexity of the response to be approximated and on the total number of unknowns in the neural network architecture, i.e., the weights and biases. The numbers of unknown weights and biases in the neural networks $\{P^{(i)}\}_{i=2}^{N_L}$ decrease significantly by using the transfer learning approach described in section \ref{sec:architecture}. As a result, we need a large amount of training samples when training $P^{(2)}$, which is feasible when assuming that low-fidelity solutions are computationally cheap to sample. As opposed to this, sampling high-fidelity solutions for training $P^{(N_L)}$ is computationally expensive, but we require significantly less training samples, due to variance decay in $\mbf{e}^{(i)}$ for increasing $i$ \cite{giles2015multilevel} and the small amount of unknown weights and biases induced by transfer learning. The quality of the training data is determined by the sample locations $\mbf{z}_k$ in \eqref{eq:trainingset}. As determining a proper size of the training/validation set is difficult \cite{goodfellow_deep_2016}, we employ a sampling strategy which iteratively adds new samples $\mbf{u}^{(i)}(\mbf{z}_k)$ when the neural network does not generalise well to the validation set.

The training procedures for the neural networks $\{P^{(i)}\}_{i=2}^{N_L}$ are explained in more detail below.
\subsubsection*{Training $P^{(2)}$}
\label{sec:trainP2}
\noindent The first neural network that needs to be trained is $P^{(2)}$. When training $P^{(2)}$, we approximate the error
\begin{linenomath*}\begin{equation}
\mbf{e}^{(2)}(\mbf{z})|_{X^{(1)}} = \mbf{u}^{(2)}(\mbf{z})|_{X^{(1)}} - \mbf{u}^{(1)}(\mbf{z})\ ,
\end{equation}\end{linenomath*}
and to construct $T^{(2)}$, the sample locations $\mbf{z}_k$ need to be specified. As the low-fidelity solutions $\mbf{u}^{(1)}(\mbf{z})$ and $\mbf{u}^{(2)}(\mbf{z})$ are assumed to be computationally cheap to sample from, we place many samples on the first two levels. In this work, the initial sample set comprises a relatively large set of $10^{\dim(I)}$ solutions in the random space $I$, which shows to give good results for our test cases in section \ref{sec:Results}. However, as we employ an iterative sampling strategy, the final size of the sample set is tailored to the underlying problem, which makes a proper size of the initial sample set less significant. The locations are defined using Monte Carlo sampling according to the underlying PDF of the uncertainties $\rho(\mbf{z})$. As transfer learning can not be used for training $P^{(2)}$, the number of unknown weights and biases is large compared to $P^{(i)},\ i=3,...,N_L$. As a result, the initial sample set should be large enough to properly train the neural network, but also small enough to be feasible in high-dimensional random spaces. The sample set is split using the 80/20\%-rule in a training set $T^{(2)}$ and a validation set $V^{(2)}$. Splitting the full sample set according to this 80/20\% ratio is known as the Pareto Principle \cite{goodfellow_deep_2016} and is commonly used in machine learning for determining the size of the training and validation sets. After minimising \eqref{eq:costfunctionlambda}, we tune the hyperparameters  and pick the network architecture which results in the smallest validation error $v_{\min}^{(i)}$. This is discussed in more detail in section III. Training is stopped if $v_{\min}^{(i)}$ satisfies:
\begin{linenomath*}\begin{equation}
v_{\min}^{(i)} < \varepsilon\ ,
\label{eq:stoppingcriterion}
\end{equation}\end{linenomath*}
where $\varepsilon$ is a specified threshold. The threshold $\varepsilon$ indicates how well the trained neural network generalises to data that is not contained in the training set. If the stopping criterion is not met, we increase the sample set by sampling another $10^{\dim(I)}$ new PDE solutions. The newly obtained larger sample set is again split randomly in a training/validation set following the 80/20\%-rule, which are used to retrain the neural network. This process of training, validating, enlarging the training set is repeated until \eqref{eq:stoppingcriterion} is satisfied. Notice that each time the sample set is increased, the neural network needs to be retrained, which is assumed to be a relatively cheap operation when compared to computing a high-fidelity solution, as the weights and biases of the previous training can be used as a very good initial guess for the retraining.

\subsubsection*{Training $P^{(i)},\ i>2$}
\label{sec:trainPi}
\noindent Constructing a large training set when increasing $i$ is often not feasible. Therefore, the combination of transfer learning and iterative sample set construction is used to limit the required number of training samples. Training the relatively small neural networks is fast in general. We iteratively increase the size of the training set during the training procedure. We again start with a small sample set that is used for training/validation of the neural network. After training the neural network we check if \eqref{eq:stoppingcriterion} is met, if not, then new samples are added to the sample set, which effectively increases the size of the training/validation set. This process is repeated until \eqref{eq:stoppingcriterion} is satisfied. 

Initially, the sample set comprises $2^{\dim(I)}$ randomly placed samples in the space $I$, which is split in a training set $T^{(i)}$ and validation set $V^{(i)}$ following the 80/20\%-rule. The size of the initial sample set is smaller when compared to the initial sample set used for training $P^{(2)}$ to alleviate computational burden when sampling the computationally expensive higher-fidelity solutions. The number of samples in the initial sample set is chosen to be small in order for the MLNN method to still be feasible for high-dimensional random spaces. As mentioned before, the use of an iterative sampling strategy makes the size of the initial sample set insignificant. Furthermore, notice that we only need to sample solutions for $\mbf{u}^{(i)}$ and $\mbf{u}^{(1)}$, as accurate approximations of $\mbf{u}^{(i-1)}$ can be constructed using \eqref{eq:approximatesolution} with $N_L = i-1$. 

Nesting of samples on different levels removes the need for sampling additional solutions for $\mbf{u}^{(i-1)}$, but introduces correlated expectations on the different levels, which is unwanted \cite{MLMC_path}. Strategies to reuse samples on different levels have recently been proposed \cite{MLMC_recycling}, and could be used to further enhance the MLNN method, but are not required to show the basic methodology that we propose.

After training the neural network, we check if criterion \eqref{eq:stoppingcriterion} is met. If the stopping criterion is not met, we increase the sample set by sampling $2^{\dim(I)}$ new PDE solutions. The newly obtained larger sample set is again split randomly in a training/validation set following the 80/20\%-rule, which is used to retrain the neural network. This process of training, validating, enlarging the training set, is repeated until we satisfy \eqref{eq:stoppingcriterion}.

\subsection{Hyperparameter Tuning}
\noindent When training the neural networks, values for the hyperparameters need to be specified. A complete list of hyperparameters in this work is shown below:
\begin{itemize}
\item Hyperparameters when training $P^{(2)}$:
\begin{itemize}
\item $\lambda$ (regularisation parameter)
\item $N_{\text{CNN}}$ (number of convolutional layers)
\item $N_{\text{FC}}$ (number of fully-connected layers)
\item $n$ (number of neurons per fully-connected layer)
\end{itemize}
\item Hyperparameters when training $P^{(i)},\ i>2$:
\begin{itemize}
\item $\lambda^{(i)}$ (regularisation parameter)
\item $n^{(i)}$ (number of neurons in the added fully-connected layer)
\end{itemize}
\end{itemize}
To tune these hyperparameters, we use a grid search, which specifies pre-defined values for each parameter, and constructs a tensor grid of all possible combinations of hyperparameter values. A neural network is then trained for each combination of hyperparameter values. To limit the total number of neural networks we have to train, we pick only 3 values for each hyperparameter, which are commonly used in deep-learning approaches:
\begin{itemize}
\item $\lambda \in \{0, 10^{-6}, 10^{-3}\}\ ,$ 
\item $N_{\text{CNN}} \in \{2, 4, 6\}\ ,$
\item $N_{\text{FC}} \in \{1, 3, 5\}\ ,$
\item $n \in \{\lfloor\frac{N^{(1)}}{2}\rfloor, N^{(1)}, 2N^{(1)}\}\ .$
\end{itemize}
This range of values for the hyperparameters has been shown to work well for many cases \cite{aggarwal_neural_2018,goodfellow_deep_2016, raissi_physics_2017, kim_deep_2018,ladicky_physicsforests:_2017,schmidhuber_deep_2015} and is therefore employed in this paper. 

\section{Complete algorithm}
\label{sec:Algorithm}
\begin{figure}[hbt!]
\centering
\includegraphics[width = \textwidth]{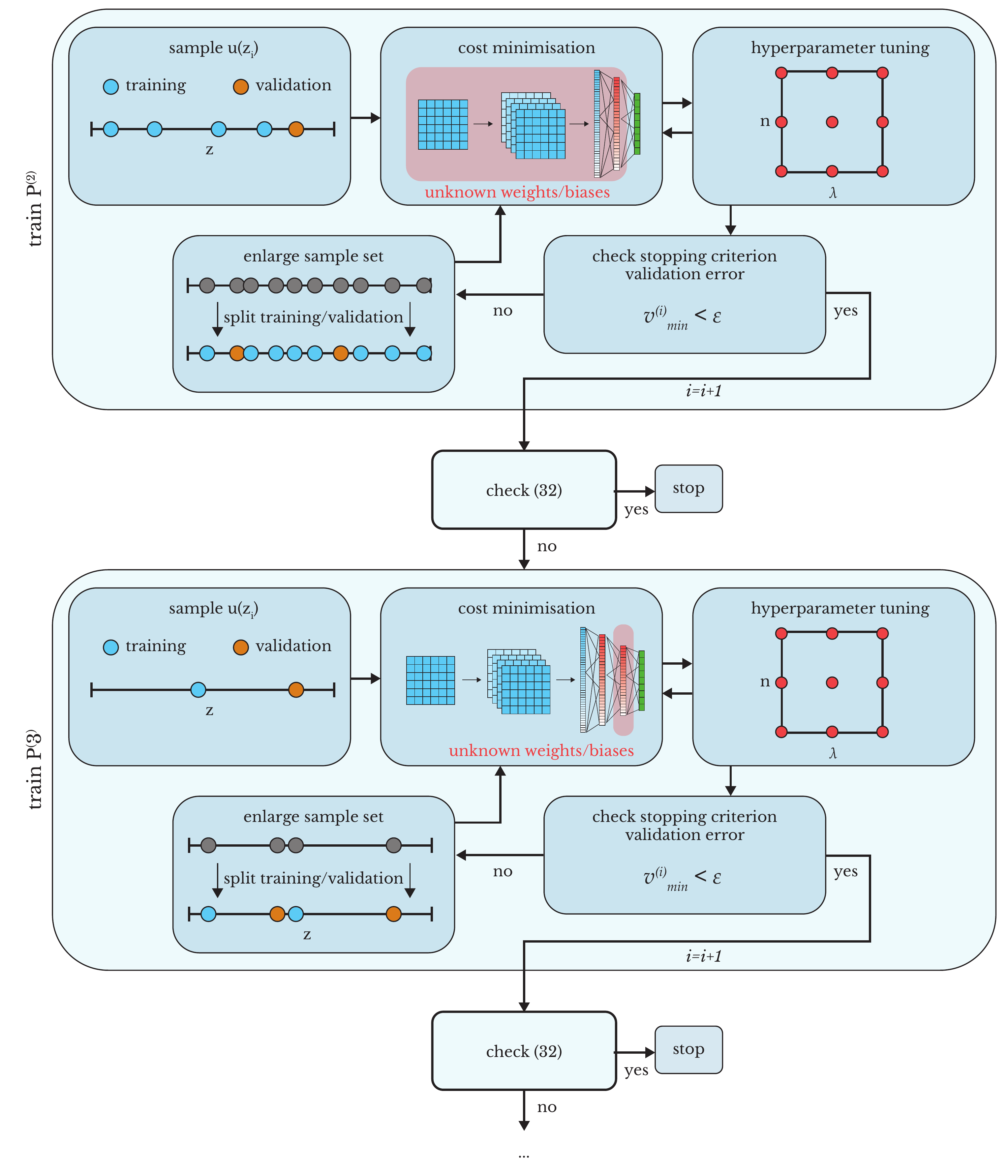}
\caption{\label{fig:completealgorithm} Schematic overview of the MLNN method.}
\end{figure}

\noindent A schematic representation of the MLNN method is shown in figure \ref{fig:completealgorithm}.
\\
The number of levels that should be picked, depends on the decay of $\mbf{e}^{(i)}$ with increasing $i$. As a result, an intuitive criterion for adding new levels is:
\begin{linenomath*}\begin{equation}
\text{if}\ \exists \mbf{u}^{(N_L-1)}\in T^{(N_L)} : \|P^{(N_L)}(\mbf{u}^{(N_L-1)})\|_2>\varepsilon_{\text{acc}} \rightarrow N_L = N_L+1\ ,
\label{eq:tolerance}
\end{equation}\end{linenomath*}
where $\varepsilon_{\text{acc}}$ is the accuracy tolerance value. To clarify, two thresholds need to be specified, i.e., $\varepsilon$ in \eqref{eq:stoppingcriterion} and $\varepsilon_{\text{acc}}$ in \eqref{eq:tolerance}. As mentioned before, the threshold $\varepsilon$ indicates how well the trained neural network generalises to data not contained in the training set. The threshold $\varepsilon_{\text{acc}}$ corresponds to the accuracy of the approximated high-fidelity solution when using \eqref{eq:approximatesolution}. In general, $\varepsilon_{\text{acc}}$ is set to the required accuracy of the approximation, and $\varepsilon$ is set to a value which is one or two orders of magnitude smaller than $\varepsilon_{\text{acc}}$ to ensure that the accuracy of the approximate high-fidelity solution is not influenced by possibly poor training of the neural networks.

The MLNN method is flexible and robust and works for a wide variety of problems, some of which are shown in the next section. We note that the assumption that $\mbf{e}^{(i-1)}(\mbf{z})$ and $\mbf{e}^{(i)}(\mbf{z})$ should possess similar behavior is a key assumption in the MLNN method. If this assumption does not hold, the transfer learning approach is not optimal and one might still require many samples for approximating $\mbf{e}^{(i)}(\mbf{z})$. The extension of transfer learning to problems where the similarity assumption does not hold will require further research in transfer learning; this is an open topic of research in the field of machine learning and it is therefore difficult to estimate how much $\mbf{e}^{(i-1)}(\mbf{z})$ and $\mbf{e}^{(i)}(\mbf{z})$ may deviate and which features are transferable \cite{Yosinski}.

\section{Results}
\label{sec:Results}
\noindent In this section we use the MLNN method to do surrogate construction for test cases with ranging complexity. These surrogates may be used to extract statistical quantities, e.g., mean and variance. As the purpose of this paper is to construct accurate surrogates, we assume that all the uncertain parameters are uniformly distributed, which does not put any emphasis on certain parameter configurations. The neural networks are constructed and trained using Tensorflow \cite{abadi_tensorflow:_2016}, which is a highly optimised library for machine learning.
\subsection{Steady-State Advection Diffusion Equation}
\noindent In this section we study the following:
\begin{itemize}
\item Construction of a parametric solution for a linear 1D differential equation.
\item Behaviour of $P^{(i)}$ for a linear differential equation.
\item Extrapolation outside the training domain.
\item Comparison with MLSC \cite{teckentrup_multilevel_2015}.
\end{itemize}
\subsubsection{Parametric Solution}
\noindent We employ the MLNN method for constructing a parametric solution for the 1D steady-state advection diffusion equation:
\begin{linenomath*}\begin{equation}
\frac{\text{d}u}{\text{d}x} - \frac{1}{\text{Re}}\frac{\text{d}^2u}{\text{d}x^2}=0,\quad u(0) = 0,\quad u(1)=1,\quad x\in[0, 1]\ ,
\label{eqAD:ADequation}
\end{equation}\end{linenomath*}
where $\text{\text{Re}}$ is the Reynolds number and is uncertain, i.e., $\mbf{z} = \text{Re}$. The exact solution is given by:
\begin{linenomath*}\begin{equation}
u_{\text{exc}}(x, \text{\text{Re}}) = \frac{\exp(x\text{Re})-1}{\exp(\text{\text{Re}})-1}\ .
\label{eqAD:exactsolution}
\end{equation}\end{linenomath*} 
We aim to construct a parametric solution for $u$ as a function of $\text{\text{Re}}\in[1, 100]$. 

\begin{figure}[!h]
\centering
\includegraphics[width =\textwidth]{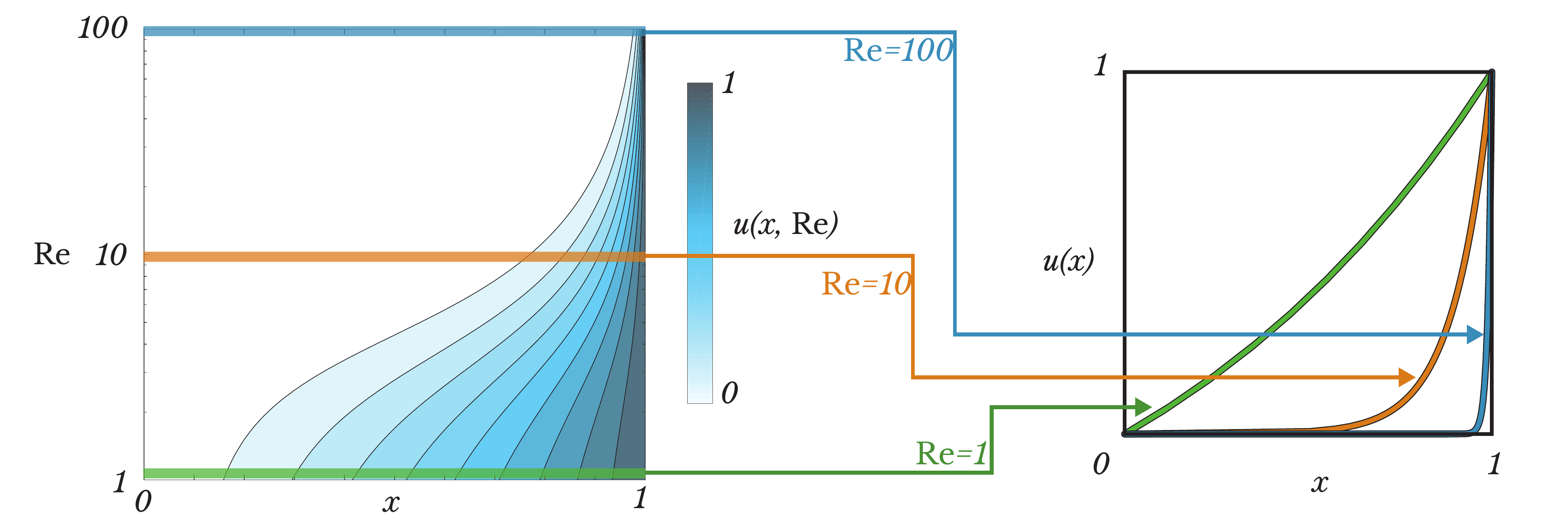}
\caption{\label{fig:SSAD_Example} Advection diffusion equation. The exact parametric solution for the advection diffusion equation.}
\end{figure}

The equations are discretised using a finite-difference approach on an equidistant grid with a resolution of $\Delta x$ with $N+1$ grid points. The solution vector on the computational grid $\mbf{u} = \left(u_i\right)_{i=0}^N \approx \left(u(x_i)\right)_{i=0}^N$, with $x_i = i\Delta x$ where $\Delta x = 1/N$, is obtained by solving the following linear system:
\begin{linenomath*}\begin{equation}
L_1\mbf{u} - \frac{1}{\text{\text{Re}}}L_2\mbf{u} = \mbf{S}(\text{\text{Re}})\ ,
\end{equation}\end{linenomath*}
where
\begin{linenomath*}\begin{subequations}\begin{align}
L_1 &= \frac{1}{2\Delta x}\left(\begin{array}{ccccc}
0&1& & & \\
-1&0&1 & & \\
 &\ddots & \ddots & \ddots & \\
 & & -1& 0 &1\\
  & & 0& -1& 0
\end{array}\right),\quad 
L_2 = \frac{1}{\Delta x^2}\left(\begin{array}{ccccc}
-2& 1& & & \\
1&-2&1 & & \\
 &\ddots & \ddots & \ddots & \\
 & & 1& -2 &1\\
  & & 0& 1& -2
\end{array}\right)\ ,\\
\mbf{S}(\text{\text{Re}}) &= \left(\begin{array}{c}
0\\
\vdots\\
\vdots\\
0\\
-\frac{1}{2\Delta x} + \frac{1}{\text{\text{Re}}\Delta x^2}
\end{array}\right)\ ,
\end{align}\end{subequations}\end{linenomath*}
where the boundary conditions enter the discretised equation via the vector $\mbf{S}$. The discretisation scheme is second order accurate and therefore the error in the solution converges with $O(\Delta x^2)$. The exact parametric solution for $\text{\text{Re}}\in[1, 100]$ is shown in figure \ref{fig:SSAD_Example}.
The fidelity of the solution increases when $N$ increases, and we choose to increase the fidelity by increasing the number of grid points with a factor 2 for consecutive levels, i.e., $N^{(i)} = 2N^{(i-1)}$. In order for the discretisation to produce stable results, the cell Reynolds number $\text{\text{Re}}_{\Delta x} = \text{\text{Re}}\Delta x$ should satisfy $\text{\text{Re}}_{\Delta x}<2$, which is satisfied by picking $\Delta x$ sufficiently small. The largest value for $\text{\text{Re}}$ is 100, and therefore picking $N^{(1)}=100$ satisfies the cell Reynolds condition.

\begin{figure}[!h]
\centering
\includegraphics[width =\textwidth]{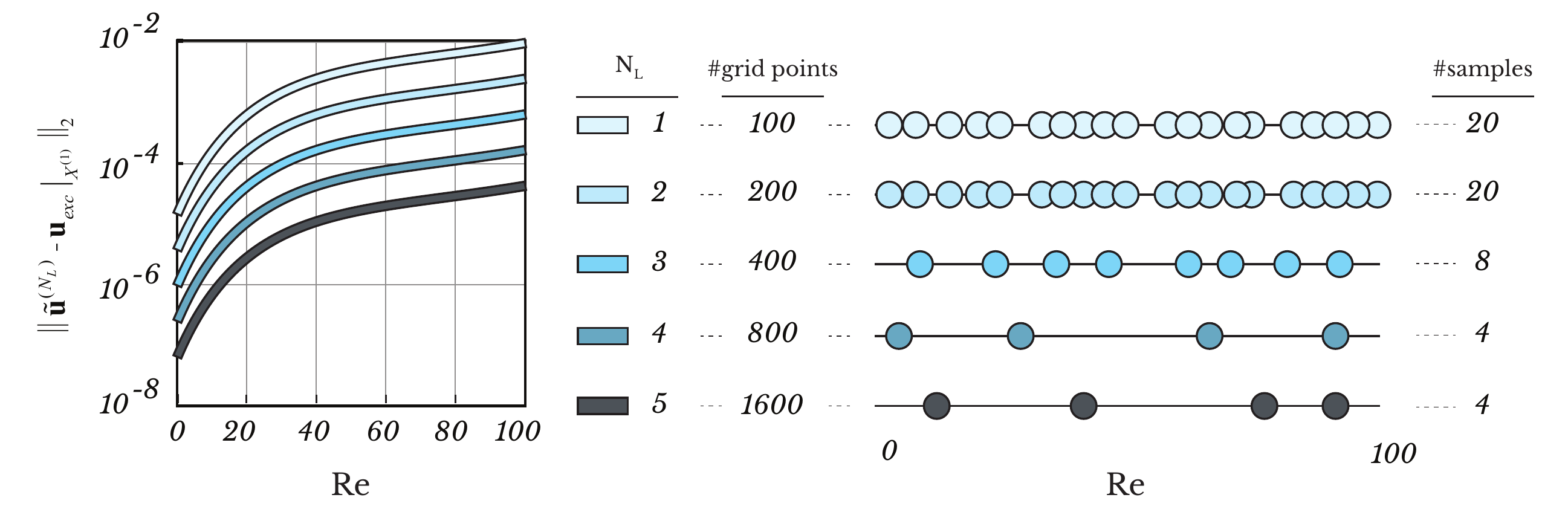}
\caption{\label{fig:SSAD_results1} Advection diffusion equation. (left) Error behaviour for different number of total levels. (right) Sample placement on each level.}
\end{figure}

The parametric solution is constructed using an accuracy tolerance $\varepsilon_{\text{acc}}=10^{-6}$, which corresponds to the desired accuracy of the surrogate. The training tolerance is set to a value which is two orders of magnitude smaller, $\varepsilon=10^{-8}$, following the rule explained in section \ref{sec:Algorithm}. As mentioned in section \ref{sec:trainP2}, we start with 10 samples on the coarsest level (8 training samples, 2 validation samples), and apply the MLNN method from there on. The error between our approximate solution \eqref{eq:approximatesolution} and the exact solution \eqref{eqAD:exactsolution} for a varying number of levels is shown in figure \ref{fig:SSAD_results1}.
The error increases with increasing Reynolds number, due to the more difficult approximation of the thin boundary layer at high Reynolds numbers. As expected, we see a decrease in required number of samples with increasing level. Furthermore, the spacing between the errors for different levels at a specific Reynolds number (left figure) indicates that the method is indeed second order accurate.

\subsubsection{Neural Network Mappings $P^{(i)}$}
\begin{figure}[!h]
\centering
\includegraphics[width =0.7\textwidth]{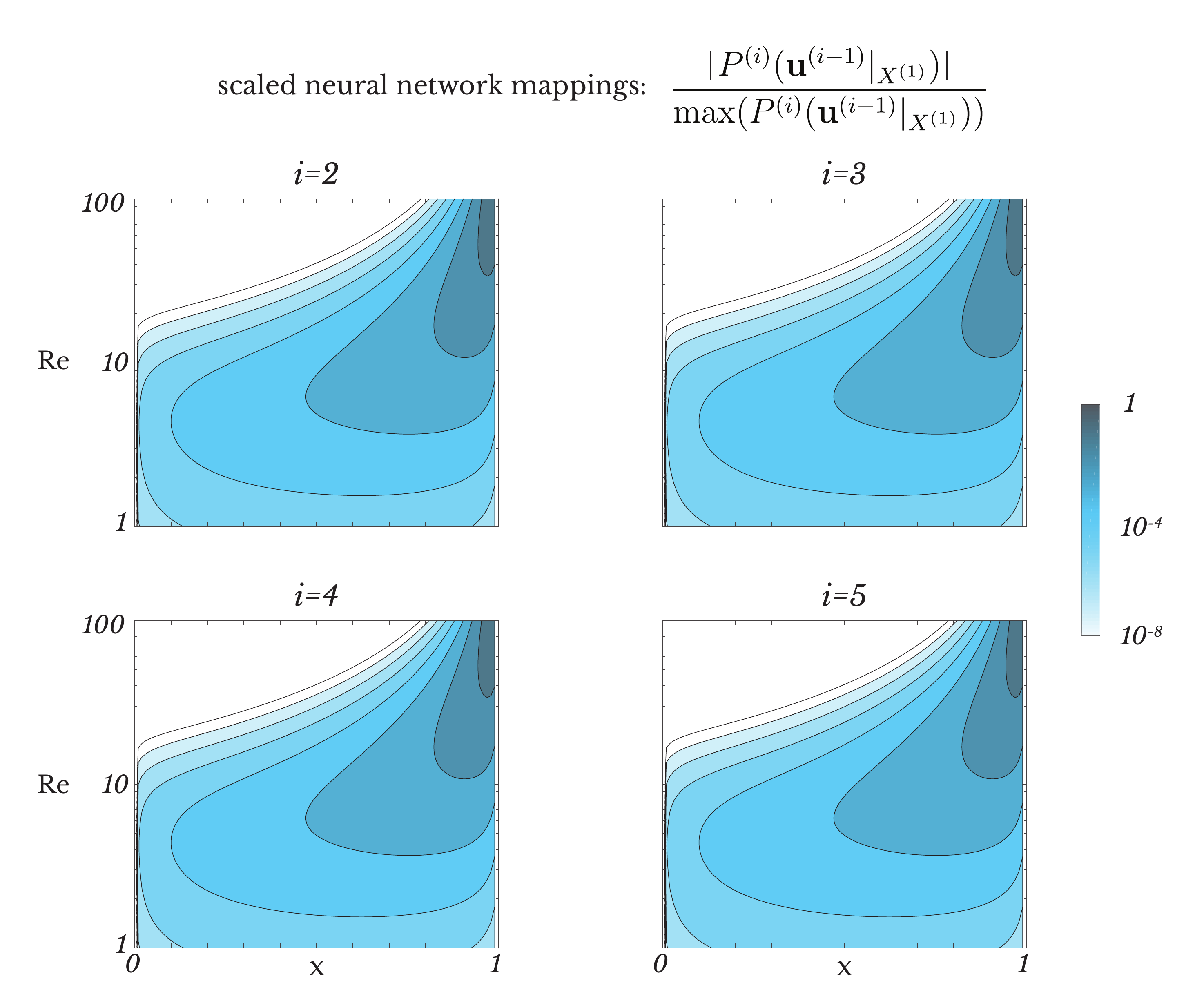}
\caption{\label{fig:SSAD_results2} Advection diffusion equation. The neural networks mappings for the advection-diffusion test case in $(x, \text{\text{Re}})$-space.}
\end{figure}

\noindent The neural network mappings $P^{(i)},\ i=2,3,4,5$ are shown in figure \ref{fig:SSAD_results2}.
Important to note is that $P^{(i)},\ i=2,3,4,5$ are approximately equal, apart from a scaling factor, which is the key insight utilised by our transfer-learning approach. This is as expected, because it can be shown that discrete solutions of this particular steady-state advection-diffusion boundary value problem satisfy the conditions outlined in \cite{proofbook} (Theorem 2.1) so that expansion (18) holds \cite{proofbook2, proofbook3}, and therefore theorem 1 holds for this case. As mentioned before, when the level increases, the error between consecutive levels decreases. As a result, less samples are required to approximate these high-level mappings, which decreases the required number of high-fidelity samples.

\subsubsection{Extrapolation Capabilities}
\begin{figure}[!h]
\centering
\includegraphics[width =\textwidth]{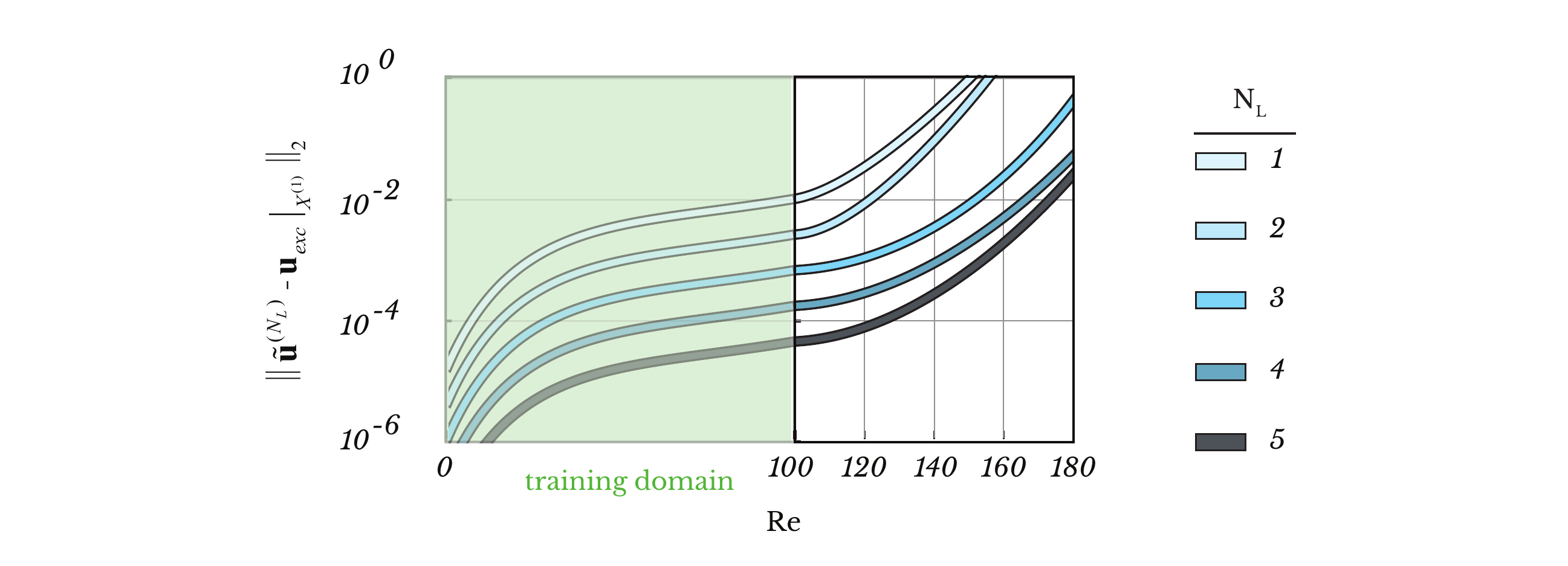}
\caption{\label{fig:SSAD_results3} Advection diffusion equation. Extrapolation errors outside the training domain $\text{\text{Re}}\in[1, 100]$.}
\end{figure}

\noindent Figure \ref{fig:SSAD_results1} shows that the MLNN method is suitable for constructing parametric solutions inside the training domain $\text{\text{Re}}\in[1, 100]$ efficiently. However, extrapolation outside the domain where the neural networks are trained is difficult in general. Figure \ref{fig:SSAD_results3} shows the errors between our approximate solution, based on $P^{(i)},\ i=2,3,4,5$ and computed with \eqref{eq:approximatesolution}, and the exact solution \eqref{eqAD:exactsolution} for $\text{\text{Re}}>100$.
The errors in the approximate solution increase rather drastically when using the trained neural networks outside the training domain. As a result, the trained neural networks are not suitable for constructing accurate solutions outside the training domain, when deviating too far from the domain boundaries \cite{goodfellow_deep_2016}. This is a common problem in machine learning. Fixing this problem is not the scope of this manuscript.

\subsubsection{Comparison with MLSC}
\begin{figure}[!h]
\centering
\includegraphics[width =\textwidth]{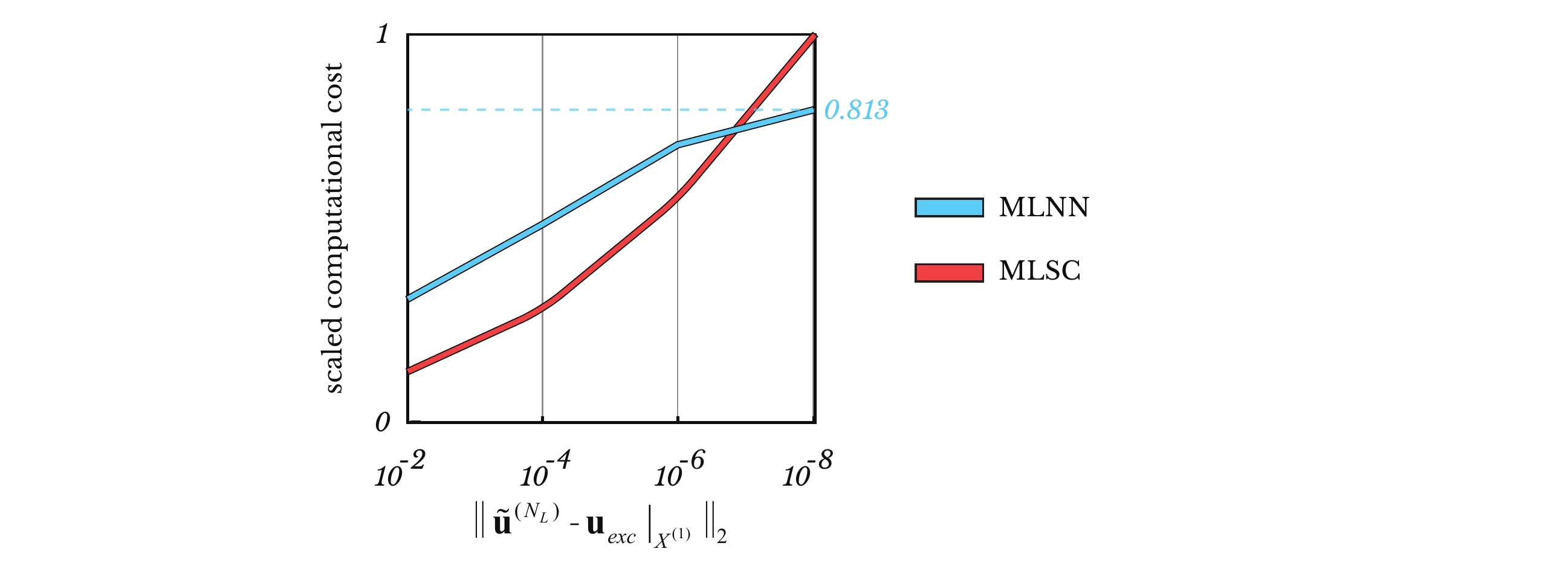}
\caption{\label{fig:SSAD_Comparison} Advection diffusion equation. The scaled computational cost when constructing a surrogate for a range of accuracies.}
\end{figure}

\noindent We compare the computational cost of MLNN and MLSC (as proposed in \cite{teckentrup_multilevel_2015}) for the case of constructing a surrogate with a specified accuracy. The computational cost comprises both solving the underlying differential equation and sampling procedure, i.e., training the neural networks (MLNN) or constructing the Clenshaw-Curtis grids (MLSC). The dependence of the computational cost on the implementation is negligible, because MLSC requires a minimum amount of implementation whereas training and constructing the neural networks for the MLNN method are performed using Tensorflow, which is a highly optimised library for machine learning. Both approaches use a sampling threshold on each level of $\varepsilon=10^{-10}$. The levels are determined by the grid resolution which is refined with a factor 2 for subsequent levels. As both approaches use the same solver with the same grid resolution on each level, we are able to compare both approaches directly. The computational costs are scaled with respect to the maximum computational cost of the MLSC approach and the results are shown in figure \ref{fig:SSAD_Comparison}.
Notice that the error on the horizontal axis is directly related to the number of levels that are used in both approaches. The MLNN method requires more computational time for the approximation on the first level when compared to the MLSC. However, the transfer learning approach significantly reduces the required number of samples on subsequent levels. As a result, MLNN appears to be more computationally efficient with increasing surrogate accuracy.

\subsection{Steady-State Burgers Equation}
\noindent In this section we increase the complexity of the previous test case by introducing a non-linearity in the underlying differential equation. We study the following:
\begin{itemize}
\item Construction of a parametric solution for a non-linear 1D differential equation.
\item Behaviour of $P^{(i)}$ for a non-linear differential equation.
\item Comparison with MLSC \cite{teckentrup_multilevel_2015}.
\end{itemize}

\subsubsection{Parametric Solution}
\noindent In order to study how the MLNN method performs in the presence of non-linearities in the underlying equation, we consider the non-linear steady-state Burgers equation:
\begin{linenomath*}\begin{equation}
\frac{1}{2}\frac{\text{d}u^2}{\text{d}x} - \frac{1}{\text{\text{Re}}}\frac{\text{d}^2u}{\text{d}x^2}=0,\quad u(0) = 0,\quad u(1)=1,\quad x\in[0, 1]\ ,
\label{eqAD:ADequation}
\end{equation}\end{linenomath*}
where $\text{\text{Re}}$ is again the Reynolds number and is assumed to be uncertain, i.e., $\mbf{z}=\text{Re}$. The effect of the non-linearity increases for increasing Reynolds number. In the linear advection-diffusion case we considered $\text{Re}\in[1, 100]$, but this range is increased to have a more pronounced non-linear effect in the solution at higher Reynolds numbers. As a result, we aim to construct a parametric solution for $u$ as a function of $\text{\text{Re}}\in[1, 1000]$.

\begin{figure}[!h]
\centering
\includegraphics[width =\textwidth]{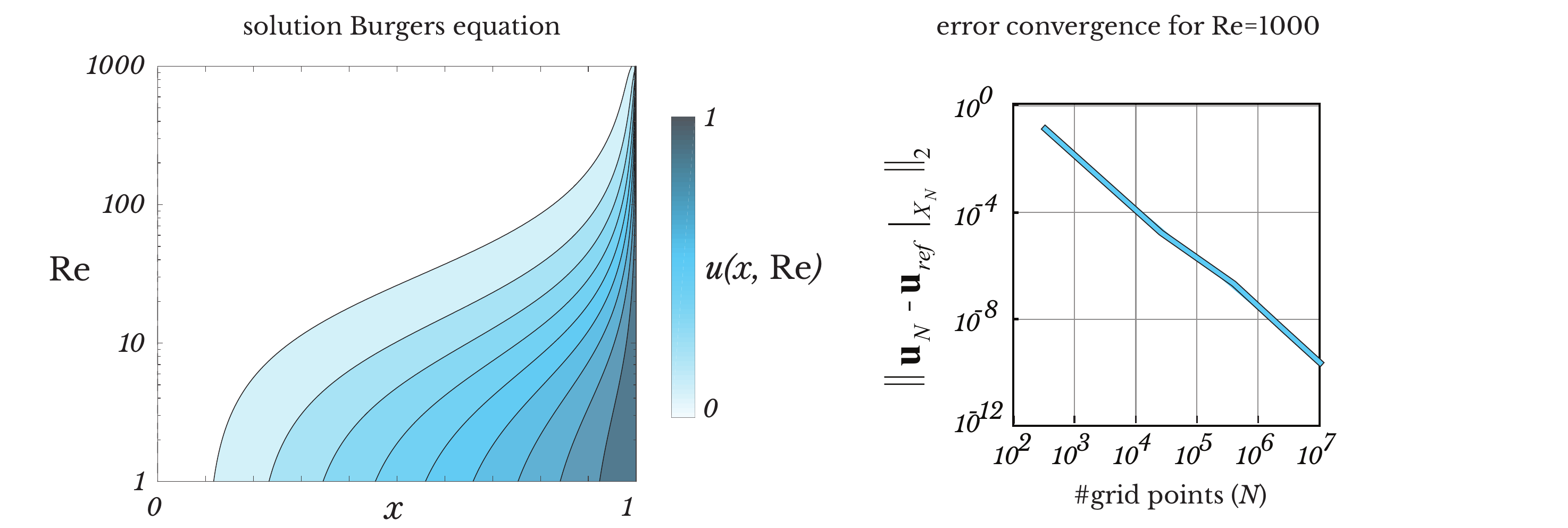}
\caption{\label{fig:SSBE_Example} Burgers equation. (left) The parametric solution for the Burgers equation computed with $N=10^5$ grid points. (right) The convergence of the error in the solution for $\text{\text{Re}}=1000$ as a function of the number of grid points. The reference solution is computed with $N=10^8$.}
\end{figure}

\begin{figure}[!h]
\centering
\includegraphics[width =\textwidth]{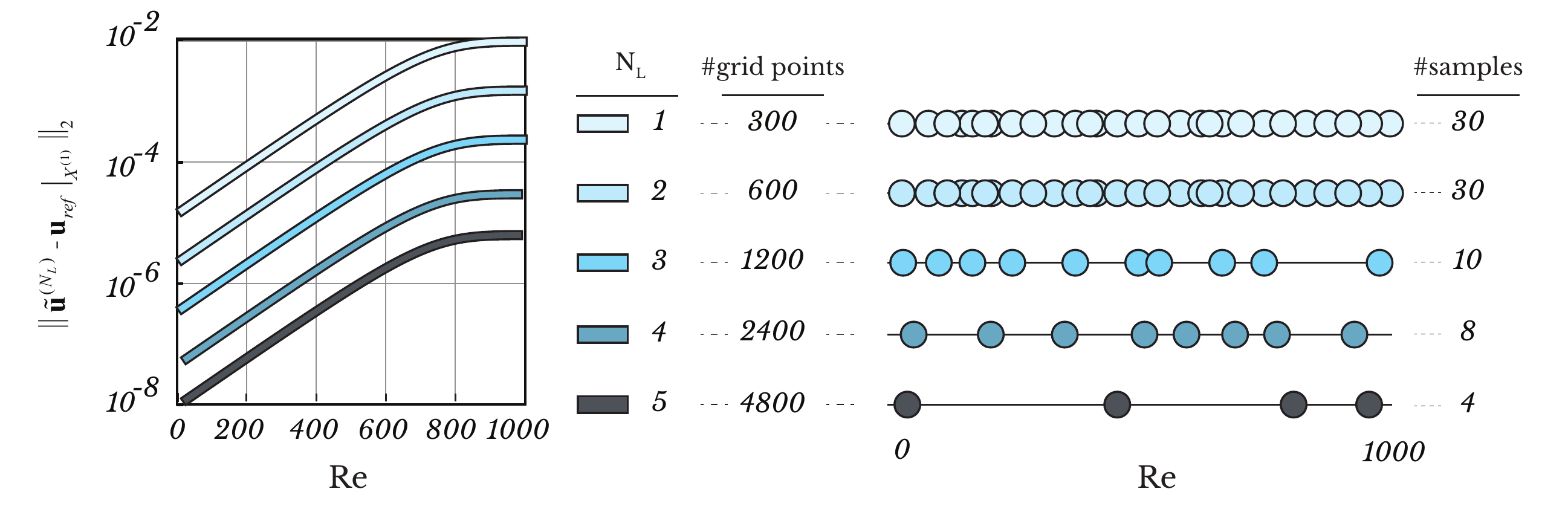}
\caption{\label{fig:SSBE_results1} Burgers equation. (left) Error convergence for different number of total levels. (right) Sample placement on each level.}
\end{figure}

The equations are discretised using a finite-difference approach on an equidistant grid with a resolution of $\Delta x$ with $N+1$ grid-points. The solution vector on the computational grid $\mbf{u} = \left(u_i\right)_{i=0}^N \approx \left(u(x_i)\right)_{i=0}^N$, with $x_i = i\Delta x$ where $\Delta x = 1/N$, is obtained by solving the following non-linear system:
\begin{linenomath*}\begin{equation}
\mbf{F}(\mbf{u}) = 0\ ,
\end{equation}\end{linenomath*}
where
\begin{linenomath*}\begin{align}
\mbf{F}(\mbf{u}) &= \left(\begin{array}{c}
u_0\\
\frac{1}{4\Delta x}(u_2^2 - u_0^2) - \frac{1}{\text{\text{Re}}\Delta x^2}(u_2 - 2u_1 + u_0)\\
\vdots\\
\frac{1}{4\Delta x}(u_N^2 - u_{N-2}^2) - \frac{1}{\text{\text{Re}}\Delta x^2}(u_N - 2u_{N-1} + u_{N-2})\\
u_N - 1
\end{array}\right)\ .
\end{align}\end{linenomath*}
This non-linear system is solved using Newton iteration and the complete discretisation scheme is second order accurate. We choose to increase the fidelity by increasing the number of grid-points with a factor 2 for consecutive levels, i.e., $N^{(i)} = 2N^{(i-1)}$, and picking $N^{(1)}=300$. The grid-resolution of the first level is picked such that we produce stable results for $\text{\text{Re}}\in[1, 1000]$. The parametric solution for $\text{\text{Re}}\in[1, 1000]$ is shown alongside the solution error convergence for $\text{\text{Re}}=1000$ in figure \ref{fig:SSBE_Example}.
As for the the linear case, the tolerance for the training procedure is set to $\varepsilon=10^{-8}$ and the accuracy tolerance is set to $\varepsilon_{\text{acc}}=10^{-6}$. As mentioned in section \ref{sec:trainP2}, we start with 10 samples on the coarsest level (8 training samples, 2 validation samples), and apply the MLNN method from there on. The error between our approximate solution \eqref{eq:approximatesolution} and the converged solution, with a varying number of levels, is shown in figure \ref{fig:SSBE_results1}. Again, we see a decrease in required number of samples with increasing level.
The error increases with increasing Reynolds number, due to the more difficult approximation of the thin boundary layer at high Reynolds numbers and possibly the increasing effect of the non-linearity in the underlying equation.

\subsubsection{Neural Network Mappings $P^{(i)}$}
\begin{figure}[!h]
\centering
\includegraphics[width =0.7\textwidth]{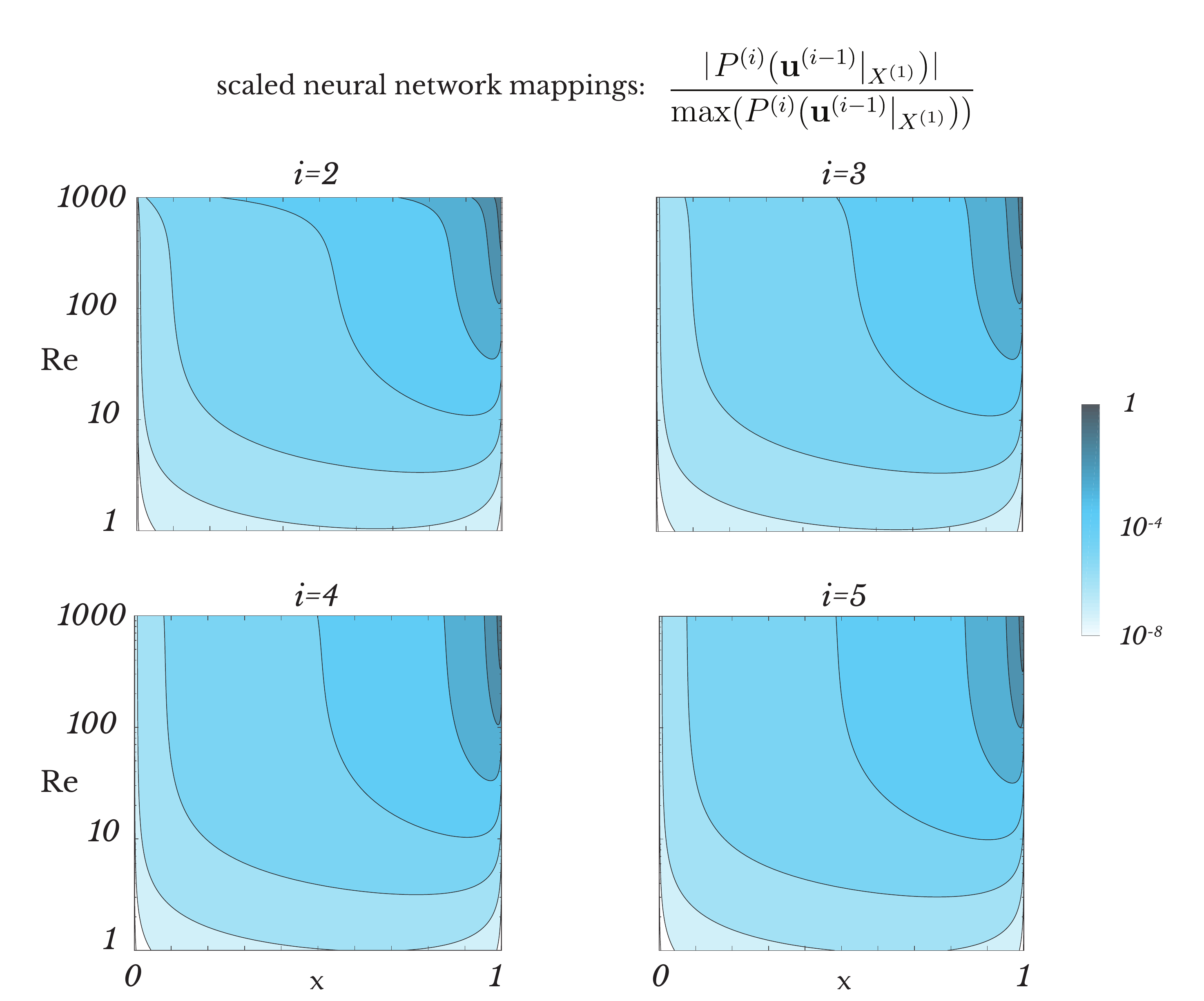}
\caption{\label{fig:SSBE_results2} Burgers equation. The neural networks mappings for the Burgers test case in $(x, \text{\text{Re}})$-space.}
\end{figure}

\noindent The neural network mappings $P^{(i)},\ i=2,3,4,5$ are shown in figure \ref{fig:SSBE_results2}.
The approximated mappings $P^{(i)}$ differ slightly for consecutive levels, which is caused by the increasing effect of the non-linearity in the underlying PDE for higher Reynolds numbers. This is in contrast to the previous linear test case, where the scaled mappings remained basically indistinguishable for increasing levels. Note that for Re $\in [1,100]$ the mappings for the consecutive levels are very similar, just as with the linear advection-diffusion equation.

To summarise, the MLNN method shows similar result when comparing it to the linear advection-diffusion test case. The number of samples required on each level increased, due to the more complex error responses $\mbf{e}^{(i)}(\mbf{z})$, which is caused by the increasing non-linearity in the underlying PDE with increasing Reynolds number. However, as the error responses are still similar, our proposed transfer learning approach is a very efficient way to learn the neural network mappings with increasing $i$.

\subsubsection{Comparison with MLSC}
\begin{figure}[!h]
\centering
\includegraphics[width =\textwidth]{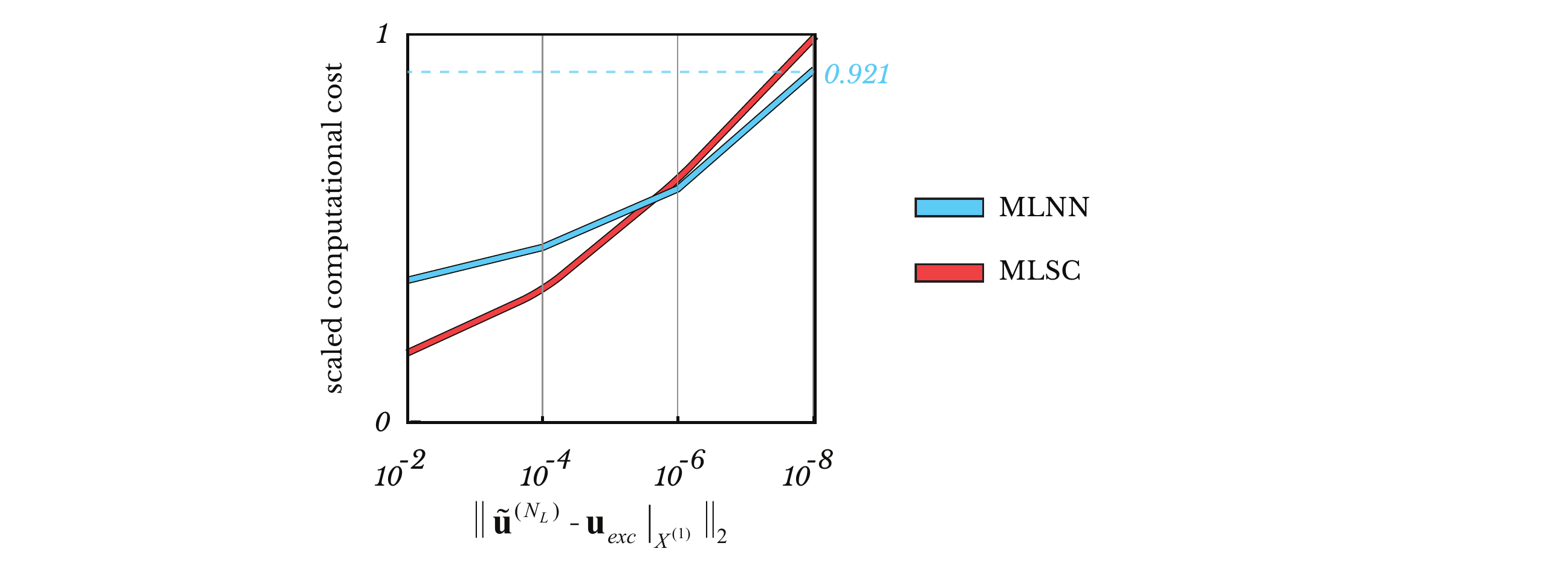}
\caption{\label{fig:SSBE_Comparison} Burgers equation. The scaled computational cost when constructing a surrogate for a range of accuracies.}
\end{figure}

\noindent We compare the computational cost of MLNN and MLSC for the case of constructing a surrogate with a certain accuracy. The computational cost is computed in the same way as described in section 6.I.4 and the results are shown in figure \ref{fig:SSBE_Comparison}.
The MLNN method requires significantly more computational time for the approximation on the first level when compared to the MLSC. Again we benefit using our approach for higher accuracy levels.

\subsection{Steady-State Incompressible Navier-Stokes Equations}
\noindent In this section we study the following:
\begin{itemize}
\item Construction of a parametric solution for the 2D steady-state incompressible Navier-Stokes equations.
\item Behaviour of $P^{(i)}$ for the steady-state incompressible Navier-Stokes equations.
\item Comparison with MLSC \cite{teckentrup_multilevel_2015}.
\end{itemize}

\subsubsection{Parametric Solution}
\noindent This test case discusses the steady-state flow over a backward-facing step, which is a common fluid mechanics problem. The governing equations are the steady-state incompressible Navier-Stokes equations in dimensionless form:
\begin{linenomath*}\begin{subequations}\begin{align}
\nabla \cdot\mbf{u} &= 0\ ,\label{eq:NSsteadycompr}\\
(\mbf{u}\cdot \nabla)\mbf{u} &= -\nabla p + \frac{1}{\text{\text{Re}}}\nabla^2 \mbf{u}\ .\label{eq:NSsteadymomentum}
\end{align}\end{subequations}\end{linenomath*}
where $\mbf{u}=(u,v)$ is the velocity field, $p$ the modified pressure, and $\text{\text{Re}}$ the Reynolds number which is assumed to be uncertain, i.e., $\mbf{z}=\text{Re}$. The incompressible Navier-Stokes equations are hard to solve due to the non-linearity and the implicit coupling between mass conservation \eqref{eq:NSsteadycompr} and momentum conservation \eqref{eq:NSsteadymomentum} by means of the pressure.

\begin{figure}[hbt!]
\includegraphics[width = \textwidth]{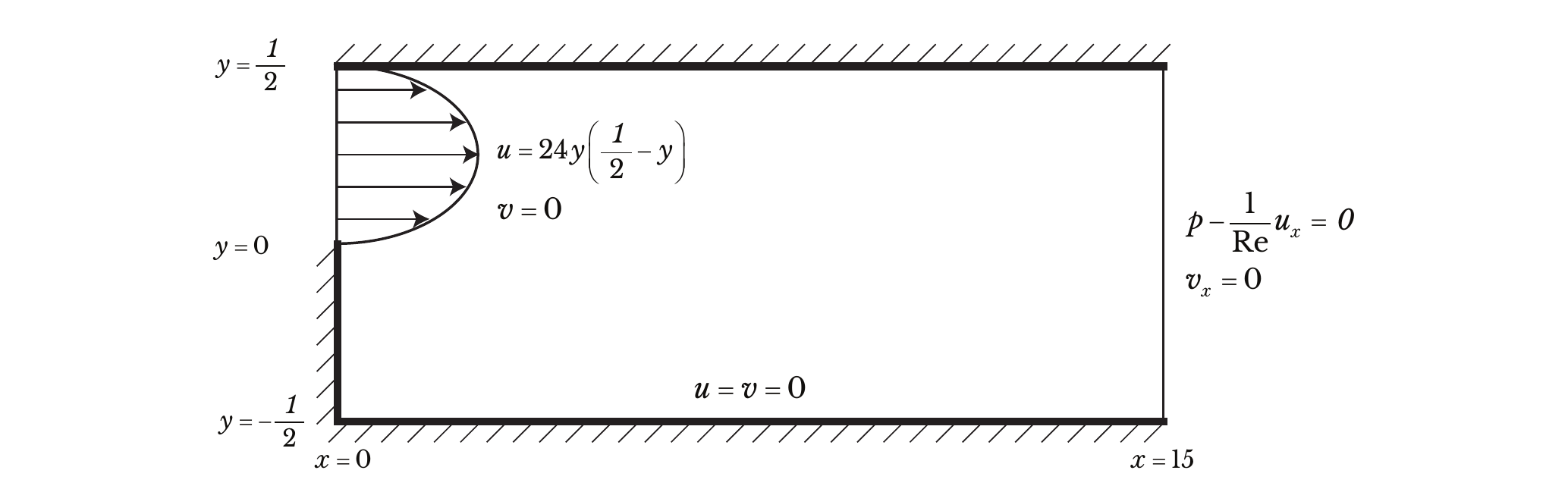}
\caption{\label{fig:BFS} Schematic of the backward-facing step problem.}
\end{figure}

\begin{figure}[hbt!]
\includegraphics[width = \textwidth]{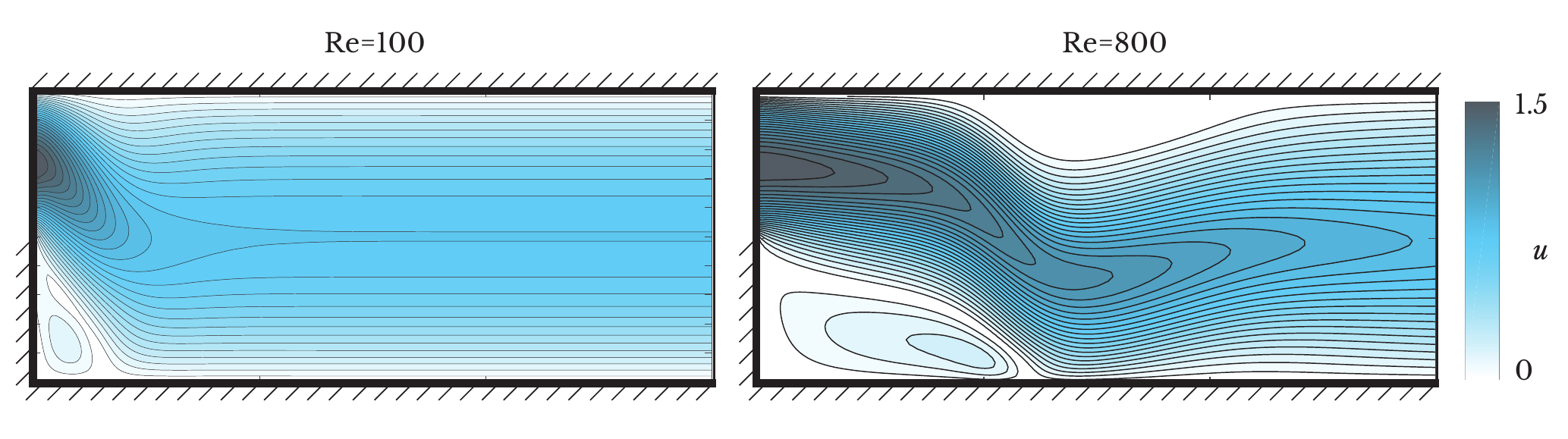}
\caption{\label{fig:BFSflows} Examples of flows over a backward facing step at different Reynolds numbers.}
\end{figure}

\begin{figure}[hbt!]
\includegraphics[width = \textwidth]{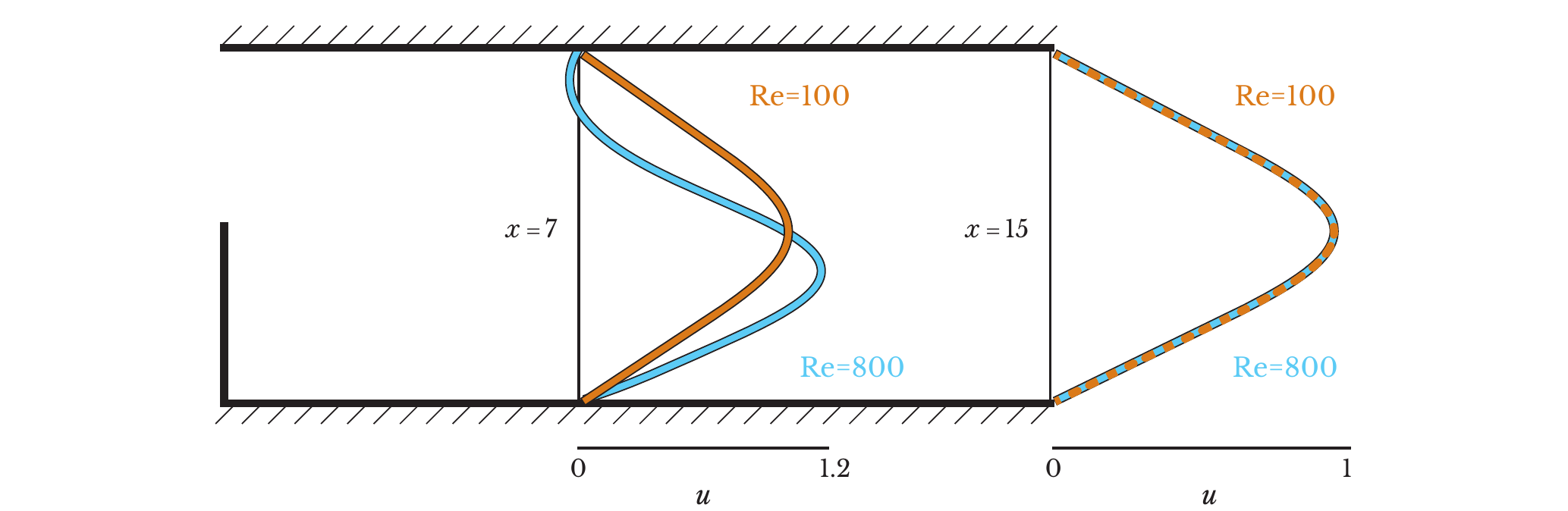}
\caption{\label{fig:BFSLFHF} The $u$-velocity profiles for backward facing step at $x=7$ and $x=15$ for two different Reynolds numbers.}
\end{figure}

The set of steady state Navier-Stokes equations \eqref{eq:NSsteadycompr}-\eqref{eq:NSsteadymomentum} are accompanied with a proper set of boundary conditions on a specified domain. In this section we consider the boundary conditions and domain that correspond to the backward-facing step problem, of which a schematic representation is shown in figure \ref{fig:BFS}.
The domain comprises a rectangle with length 15 and height 1 with solid boundaries, except for the upper part of the inlet ($x=0$) and the full outlet ($x=15$). A Poiseuille flow is imposed at the upper part of the inlet, while a pressure outlet condition is enforced at the outlet. The solver \cite{benjaminraport} is verified by comparing solutions for different mesh-sizes with a reference solution, which is believed to be an accurate benchmark for \text{Re}$=800$ \cite{Gartling}. The flow is highly dependent on the Reynolds number, and the solutions for $\text{\text{Re}}=100$ and $800$ are shown in figure \ref{fig:BFSflows}.
The Reynolds number determines the size and location of the recirculating flow region right after the step and near the top boundary. The flow develops to a steady Poiseuille flow after a distance which is determined by the Reynolds number. However, in this case the flow is not yet fully developed and that is why we impose no Poisseuille boundary conditions at the outlet. The $u$-velocity profiles at $x=7$ and $x=15$ are shown in figure \ref{fig:BFSLFHF}.

\begin{figure}[!h]
\centering
\includegraphics[width =\textwidth]{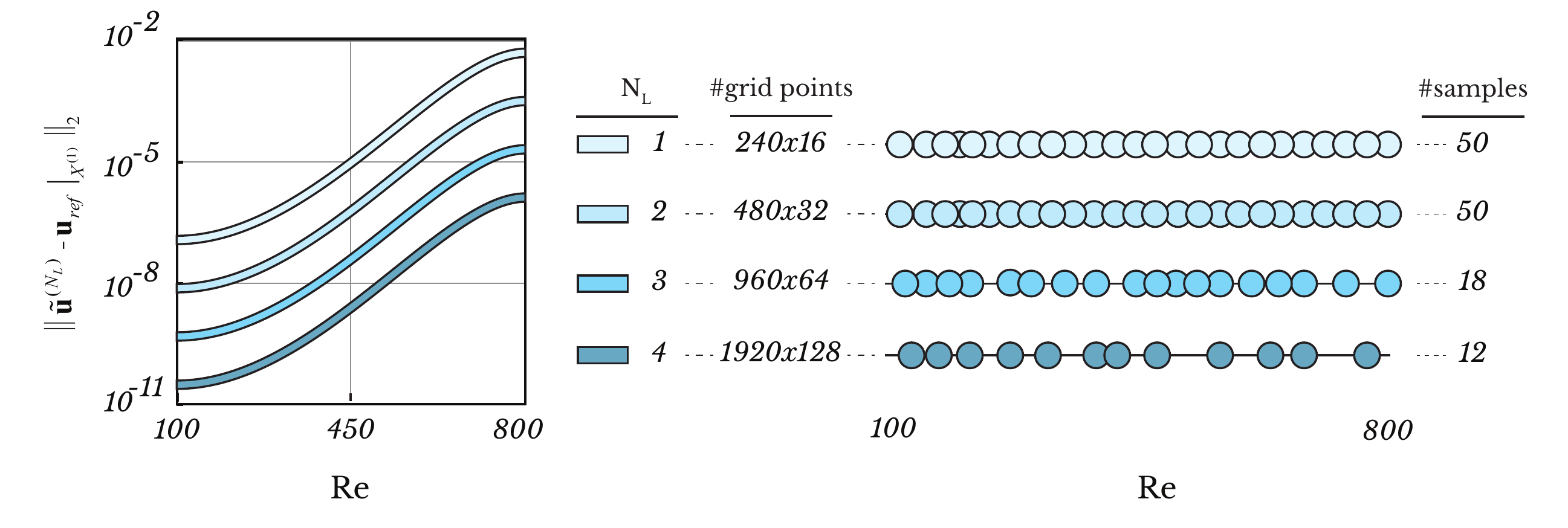}
\caption{\label{fig:SSNS_results1} Navier-Stokes equations. (left) Error convergence for different number of total levels. (right) Sample placement on each level.}
\end{figure}

The solutions on the first level are computed on a $240\times 16$ grid, and consist of the three quantities $(\mbf{u}, \mbf{v}, \mbf{p})$. We choose to increase the fidelity by increasing the number of grid points with a factor 2 for consecutive levels, i.e., $(N_x^{(i)}, N_y^{(i)}) = 2(N_x^{(i-1)}, N_y^{(i-1)})$. All three quantities $(\mbf{u}, \mbf{v}, \mbf{p})$ are given as input to the neural network as a 2D 3-channel convolutional input. As in the previous Burgers test case, the tolerance for the training procedure is set to $\varepsilon=10^{-6}$ and an accuracy tolerance $\varepsilon_{\text{acc}}=10^{-4}$ is used. We start with 10 samples on the coarsest level (8 training samples, 2 validation samples), and apply the MLNN method from there on. The normed error difference for different numbers of levels, between our approximate solution and the solution computed on a fine $3840\times 256$ grid (corresponding to level 5), is shown in figure \ref{fig:SSNS_results1}.
The number of samples required for approximating the mappings accurately, increases when compared to the previous test cases. This is due to the increase in degrees of freedom in the neural networks, as the 2D multi-channel inputs require 2D convolutional layers, instead of the 1D convolutional layers that were used in the previous test cases. Furthermore, the mappings are more complex when compared to the previous two test cases.

\subsubsection{Neural Network Mappings $P^{(i)}$}
\begin{figure}[!h]
\centering
\includegraphics[width =\textwidth]{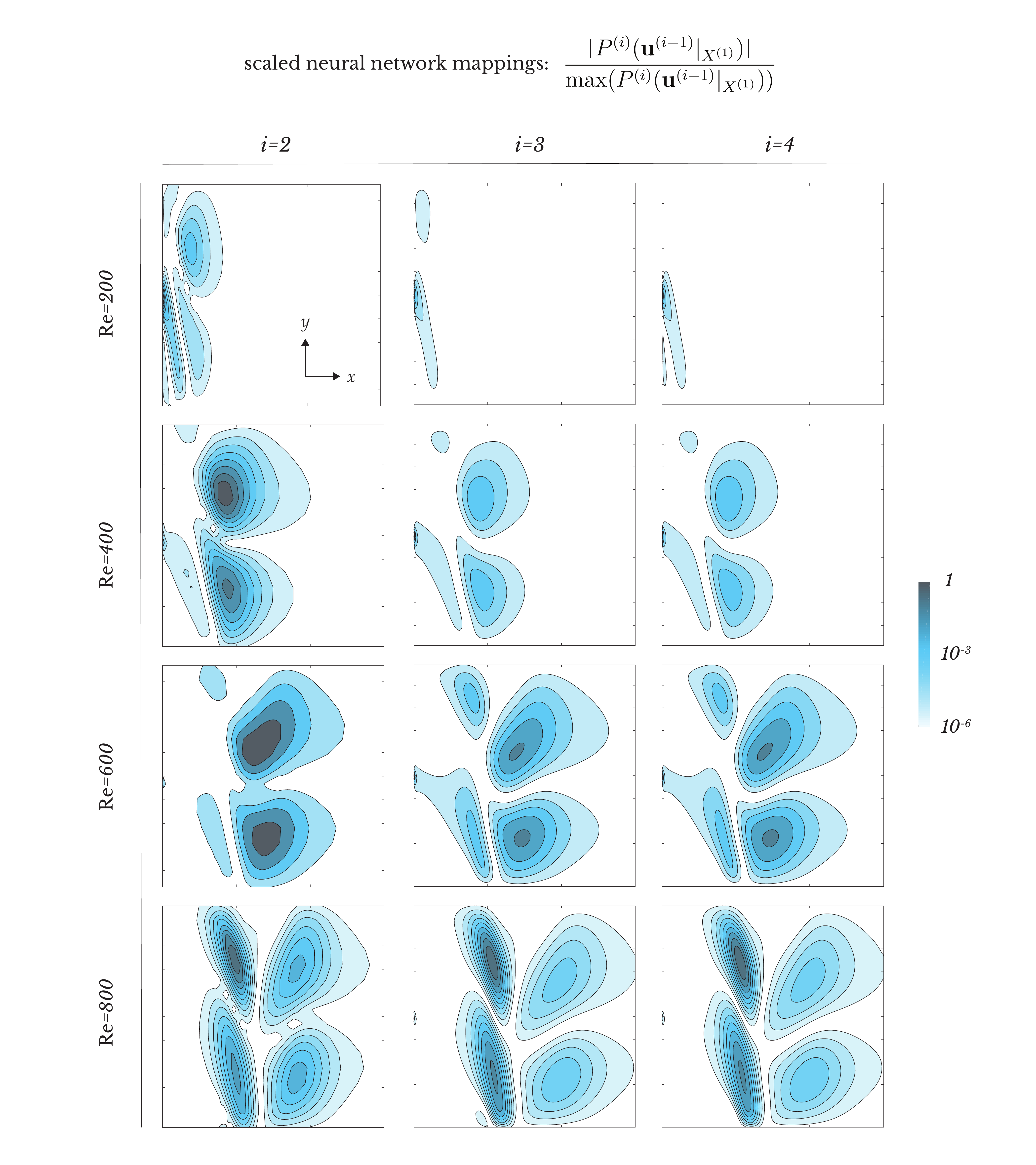}
\caption{\label{fig:SSNS_results3} Navier-Stokes equations. The neural networks mappings for the Navier-Stokes test case in $(x, y, \text{Re})$-space.}
\end{figure}

\noindent The neural network mappings $P^{(i)},\ i=2,3,4$ are shown in figure \ref{fig:SSNS_results3}.
It is striking that, even for this complex non-linear test case, the neural network approximations still show very similar behaviour for consecutive levels. We stress once again that this is the property that is utilised by our proposed transfer-learning approach to reduce the total number of samples required.

\subsubsection{Comparison with MLSC}
\begin{figure}[!h]
\centering
\includegraphics[width =\textwidth]{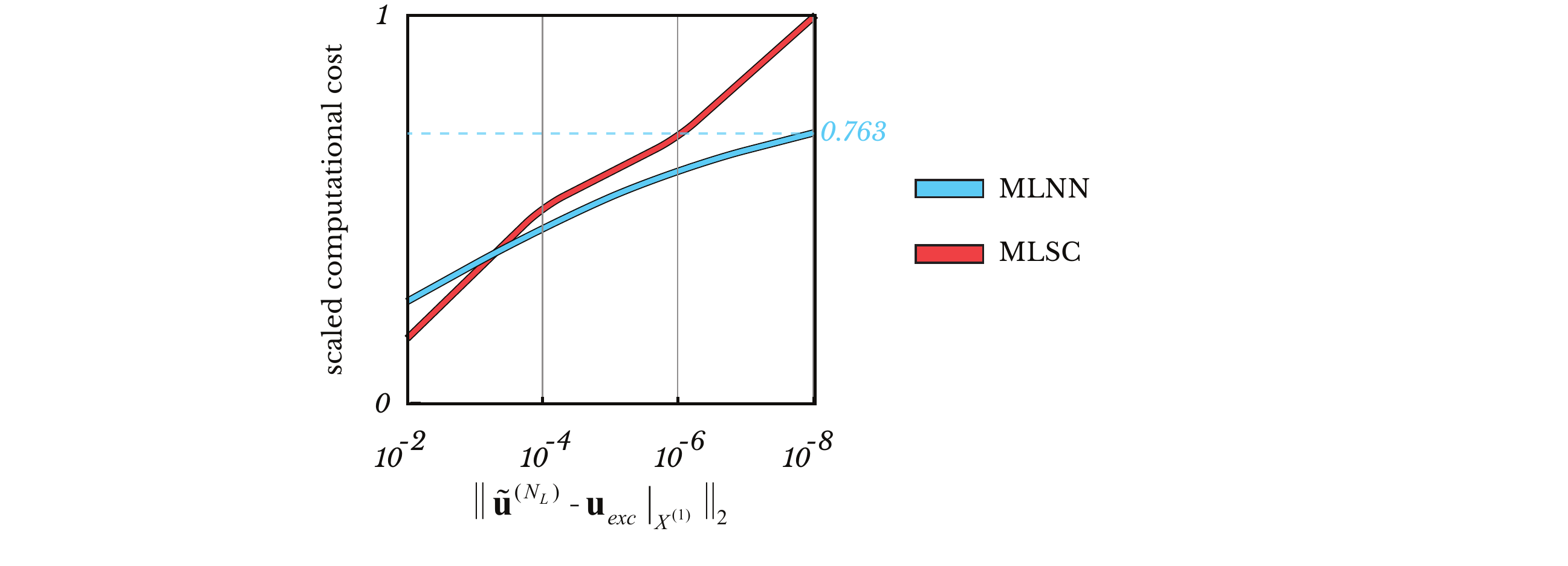}
\caption{\label{fig:SSNS_Comparison} Navier-Stokes equations. The scaled computational cost when constructing a surrogate for a range of accuracies.}
\end{figure}

\noindent We compare the computational cost for MLNN and MLSC when constructing a surrogate with a certain accuracy. The computational cost is computed in the same way as described in section 6.I.4 and the results are shown in figure \ref{fig:SSNS_Comparison}.
Again, the MLNN method requires more computational time for the approximation on the first level when compared to MLSC. However, for this non-linear test case, the MLNN method shows an advantage over the MLSC approach, which is more pronounced than for the Burgers case.

\subsection{Surrogate Modelling for 3D Fluid Sloshing}
\noindent In this section we investigate the MLNN method to construct a surrogate for unsteady fluid sloshing in a 3D rectangular tank with 2 uncertainties. These surrogates are perfectly suited for accurate uncertainty propagation. The governing equations are the 3D unsteady incompressible Navier-Stokes equations:
\begin{linenomath*}\begin{subequations}\begin{align}
\nabla\cdot \mbf{u} &= 0 \label{eqNS:incompressible}\ ,\\
\der{\mbf{u}}{t} + (\mbf{u}\cdot\nabla)\mbf{u} &= -\frac{\nabla p}{\rho} + \mbf{g} + \nu \nabla^2\mbf{u}\label{eqNS:momentum}\ ,
\end{align}\end{subequations}\end{linenomath*}
where $\mbf{u}$ is again the velocity field, $p$ the pressure, $\rho$ the density, $\mbf{g}$ the gravitational body force, and $\nu$ the kinematic viscosity. In contrast to previous test cases, the solver is not implemented in dimensionless form as we followed the implementation described in \cite{bridson_fluid_2015}.

\subsubsection*{The Solver}
\noindent In this work we use a PIC/FLIP solver \cite{zhu_animating_2005, bridson_fluid_2015} for predicting single-phase free-surface flows. PIC/FLIP is uses elements of a particle-based method such as Smoothed Particle Hydrodynamics (SPH) \cite{crespo_dualsphysics:_2015} (kernel approximation of velocities) and a grid-based finite-volume method \cite{hirsch_numerical_2007}. The main difference between pure particle-based methods such as SPH is that the particles in PIC/FLIP are passive and interactions take place on the grid. The positions of the particles that represent the fluid are evolved over time by using velocity values that are computed on a staggered grid.

\subsubsection*{Solver Initialisation}
\begin{figure}[!h]
\centering
\includegraphics[width =\textwidth]{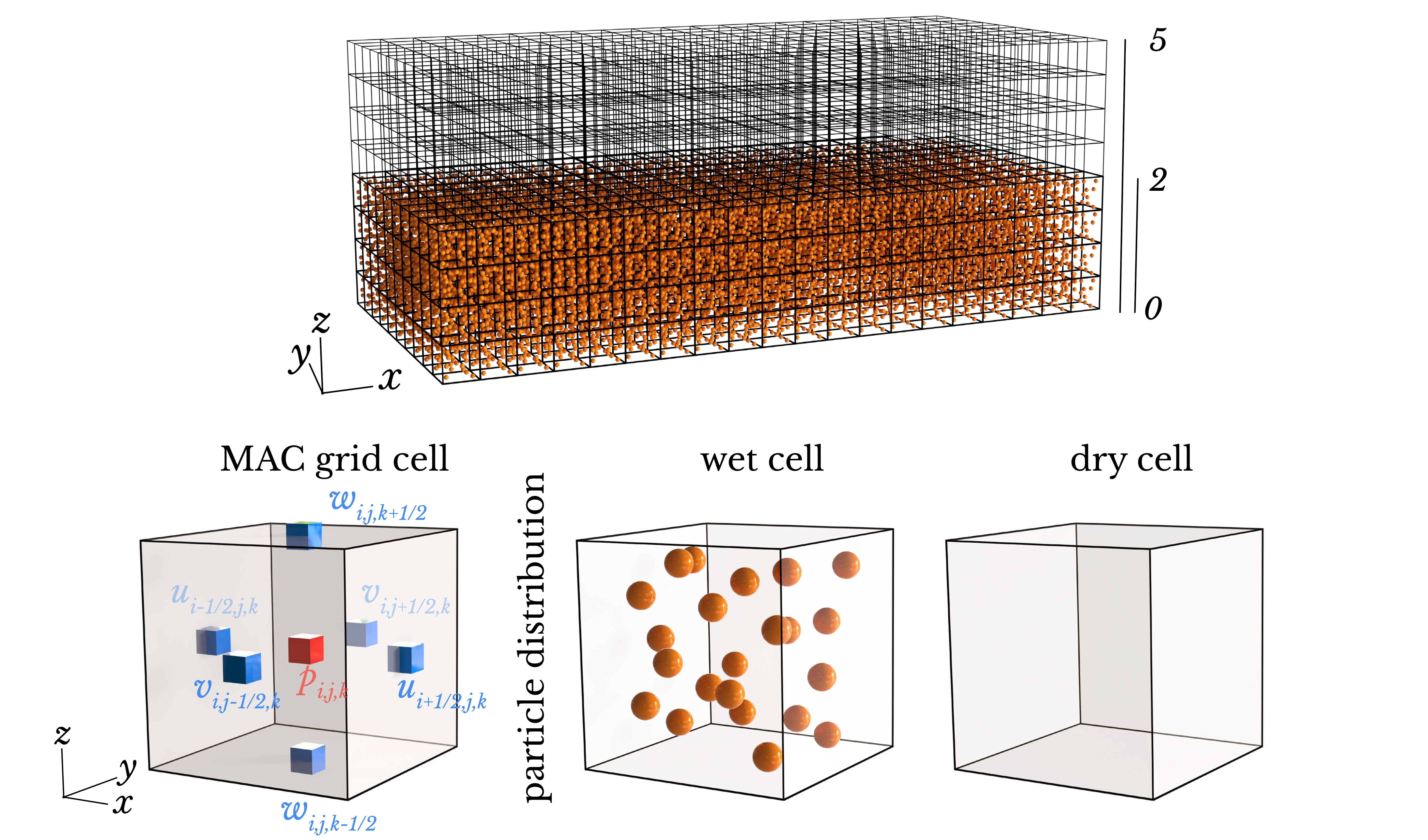}
\caption{\label{fig:USNS_init} Example of solver initialisation. A 3D staggered (128, 64, 64) grid, which is partially filled with particles.}
\end{figure}

\noindent In our simulations we consider a rectangular shaped tank of dimensions $\mbf{x}=(x,y,z)\in[0,10]\times[0,5]\times[0,5]$. The computational cost of a PIC/FLIP solver is mainly due to the operations on the grid, as the operations performed for the particles are highly parallelisable. As a result, we can use a large amount of particles $N_p=10^6$ in all simulations and randomly place these particles in the rectangular tank below $z=2$ according to a uniform distribution. The grid resolution differs depending on the fidelity of the simulation, which is discussed hereafter. An example of a particle/grid initialisation is shown in figure \ref{fig:USNS_init}.

\subsubsection*{Uncertainties}
\noindent To excite sloshing, the fluid in the partially filled tank needs forcing. In this test case we excite the fluid motion by rotating the gravitational body force. The gravity vector is defined as\begin{linenomath*}\begin{equation}
\mbf{g} = 9.81(\sin(\psi)\cos(\phi),\sin(\phi), -\cos(\psi)\cos(\phi))\ ,
\end{equation}\end{linenomath*}
where the two angles $\psi$ and $\phi$ represent rotation in the $x,y$-plane and $y,z$-plane, respectively. These angles are assumed to change during the simulation as a function of time, given by
\begin{linenomath*}\begin{equation}
\psi(t) = \left\lbrace \begin{array}{l l}
A\sin(2\lambda t),&\ t<\frac{2\pi}{\lambda}\ ,\\
0,&\ t\geq \frac{2\pi}{\lambda}\ ,
\end{array}\right.\ \ \phi(t) = \left\lbrace \begin{array}{l l}
A\sin(\lambda t),&\ t<\frac{2\pi}{\lambda}\ ,\\
0,&\ t\geq \frac{2\pi}{\lambda}\ ,
\end{array}\right.
\label{eq:motion}
\end{equation}\end{linenomath*}
where the amplitude $A$ and period $\lambda$ of the oscillation are assumed to be uncertain, i.e., $\mbf{z}=(A, \lambda)$ and are assumed to be uniformly distributed on the intervals $I_{A}=[\frac{\pi}{8}, \frac{3\pi}{8}]$ and $I_{\lambda}=[\frac{1}{2}, \frac{3}{2}]$, respectively. This particular motion leads to heavy sloshing inside the tank and is suitable for testing the applicability of the MLNN method for a highly irregular fluid motion.

\subsubsection*{The QoI}
\noindent When performing PIC/FLIP simulations, the shape of the fluid surface is our main interest, which follows directly from the particle positions. As a result, the QoI is defined as the number of particles contained within each grid cell at a given time level, which is chosen to be $t=3\pi$. At this time instant, the tank has the same orientation for all considered motion parameters, which allows us to compare the QoI.

\subsubsection*{Fidelities}
\begin{figure}[!h]
\centering
\includegraphics[width =\textwidth]{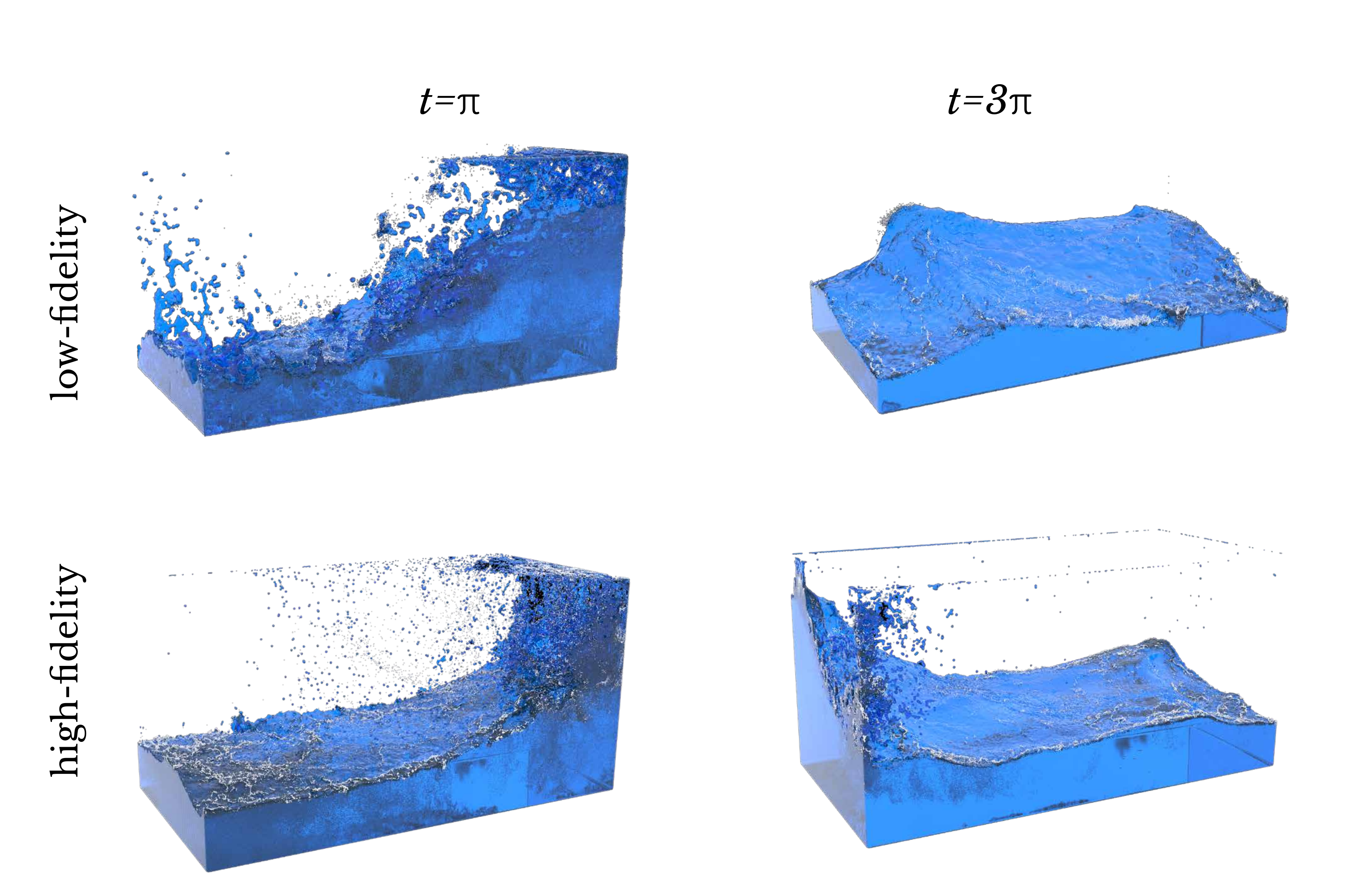}
\caption{\label{fig:SloshingExample} 3D Sloshing. Example simulation for $A=\frac{\pi}{4}$ and $\lambda=1$. Isosurface generation is used to convert the particle positions to a fluid surface, which can then be rendered \cite{crespo_dualsphysics:_2015}.}
\end{figure}

\noindent As most of the computational expense in a PIC/FLIP simulation comes from the grid based computations, we define the fidelity as the grid resolution. The time-step is tuned accordingly to satisfy the stability condition with safety factor 0.8 \cite{bridson_fluid_2015}. The solutions on the first level are computed on a fixed grid of $32\times 16 \times 16$ cells. We choose to increase the fidelity by increasing the number of grid cells in every spatial direction with a factor 2 for consecutive levels, i.e., $(N_x^{(i)}, N_y^{(i)}, N_z^{(i)}) = 2(N_x^{(i-1)}, N_y^{(i-1)}, N_z^{(i-1)})$. Obtaining a fully converged solution for this test case is difficult and we therefore consider as reference solution a solution that is computed on a fine $256 \times 128 \times 128$ grid, which corresponds to the 4th level. The QoI defined on the low-fidelity grid is given as input to the neural network as a 3D 1-channel convolutional input. An example of a low and a high-fidelity simulation result for $A=\frac{\pi}{4}$ and $\lambda=1$ is shown in figure \ref{fig:SloshingExample}.

\subsubsection*{Parametric Solution}
\begin{figure}[!h]
\centering
\includegraphics[width =\textwidth]{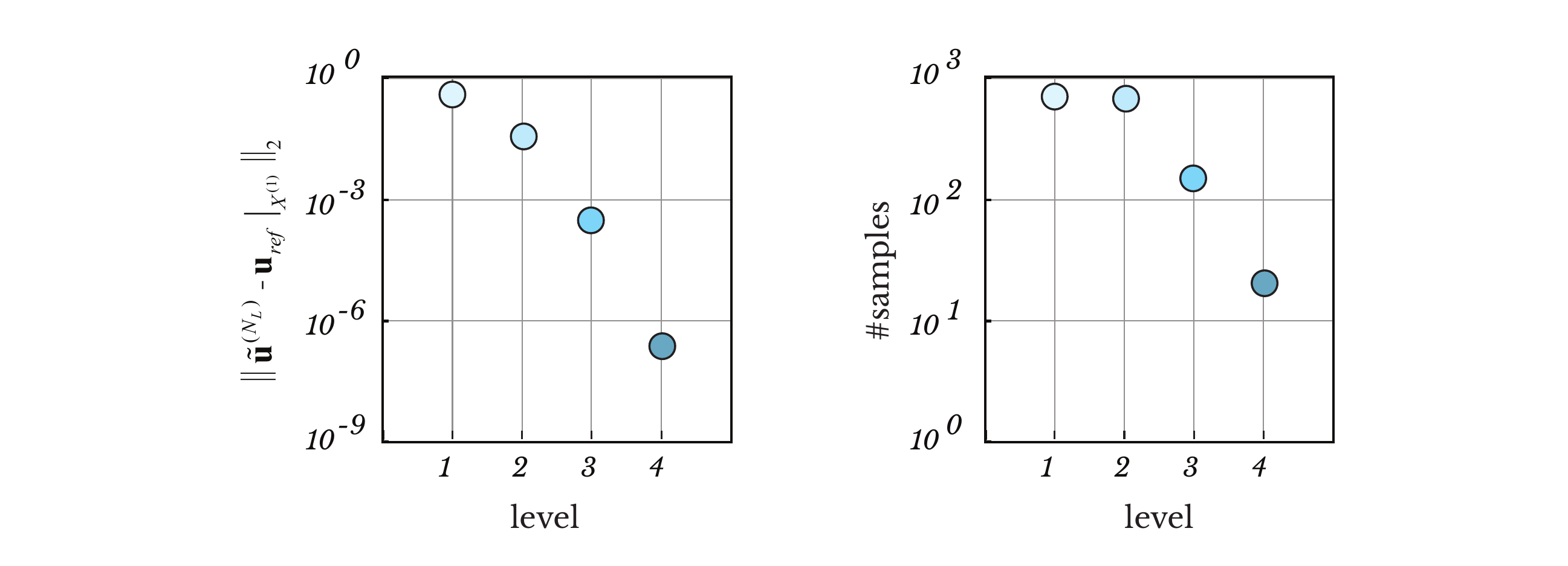}
\caption{\label{fig:SloshingResults} 3D Sloshing. Error convergence of the surrogate for the QoI constructed with the MLNN method.}
\end{figure}

\noindent The number of particles in each of the low-fidelity grid cells, i.e., the QoI, is a function of the uncertain parameters $A$ and $\lambda$. The MLNN method is used to construct a surrogate model of the QoI in the random space spanned by the uncertain parameters. The goal is to approximate the solution that is computed on the high-fidelity $256 \times 128 \times 128$ reference grid, which corresponds to the 4th level in the MLNN approach. Minimising the samples on the 4th level required for accurate surrogate construction is paramount for the feasibility of the MLNN approach. As before, the tolerance for the training procedure is set to $\varepsilon=10^{-6}$ and an accuracy tolerance $\varepsilon_{\text{acc}}=10^{-4}$ is used. Errors are computed using a reference surrogate that is constructed using high-fidelity simulations (on level 4) with a $25\times 25$ Gauss-Legendre grid \cite{xiu_numerical_2010} and the $L_2$-norm is taken over the solutions in the entire parameter space. The resulting error convergence for the QoI is shown in figure \ref{fig:SloshingResults}.
\\
For these highly non-linear fluid motions, the surrogate construction is challenging, which can be noticed in the required number of samples on each level. The figure shows that constructing the surrogate requires significantly more samples on each level when compared to previous test cases, which is caused by having two uncertainties and the complexity of the test case. However, there is still a significant decrease in the required number of samples when increasing the level, which shows the applicability of our transfer learning approach also for this complex test case.

\subsubsection*{Comparison with MLSC}
\begin{figure}[!h]
\centering
\includegraphics[width =\textwidth]{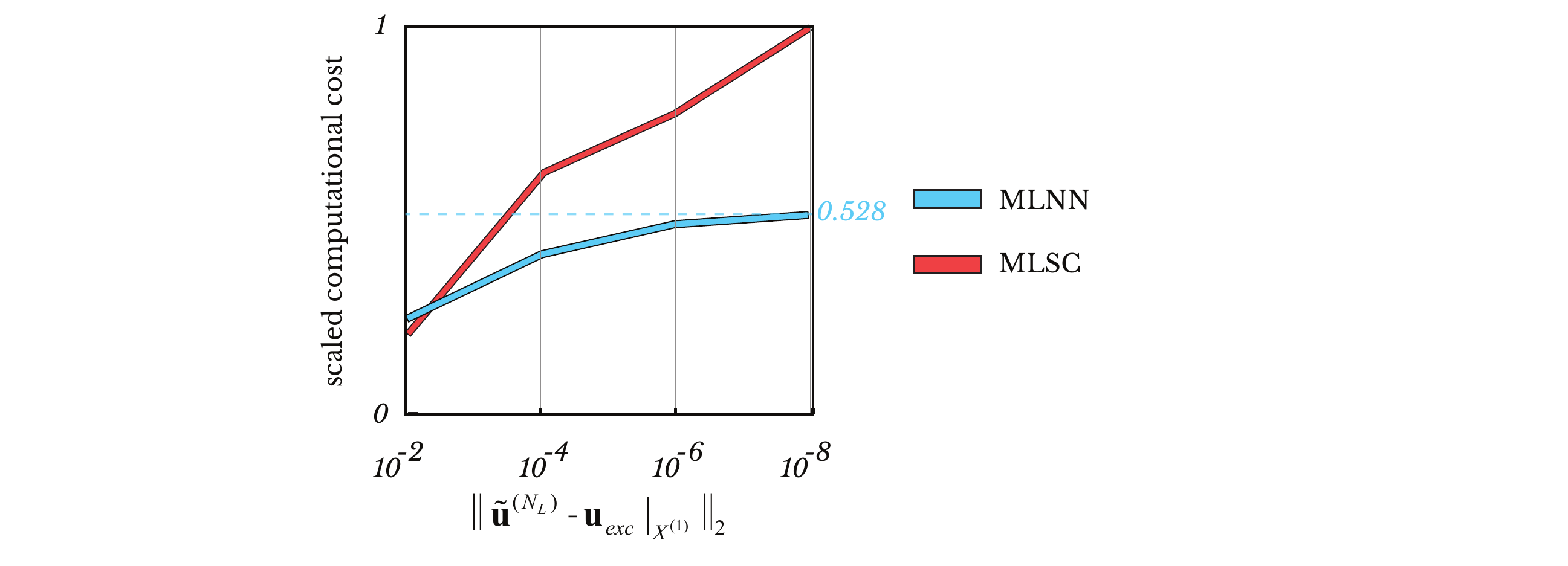}
\caption{\label{fig:SloshingResults1} 3D Sloshing. Comparison with MLSC.}
\end{figure}

\noindent The MLNN method is again compared with the MLSC method and the results are summarised in figure \ref{fig:SloshingResults1}.
\\
Our approach requires slightly more computational time for the approximation on the first level when compared to MLSC. On second and higher levels, our method shows a significant gain in computational efficiency when compared to the MLSC method, which is more pronounced than in previous test cases. In this case, the MLNN method has a clear advantage when used for surrogate construction with medium to high level of accuracy.

\section{Conclusion}
\noindent In this paper we have presented a novel method for surrogate construction. The method is based on a multi-level expansion of the solution and constructs approximations of the relative global discretisation error on different levels to enhance low-fidelity solutions. Inspired by the idea that these errors can be expressed in terms of the solution, the approximations of the error on each level are constructed using convolutional neural networks that apply a non-linear mapping of the solution values to the difference between the solutions computed on the corresponding level and the subsequent level. Transfer learning reduces the amount of training samples that are required to properly train the neural networks.

The MLNN method has been employed for surrogate construction for several parametric partial differential equations: steady 1D advection diffusion equation, steady 1D Burgers equation, steady 2D incompressible Navier-Stokes equations, and unsteady 3D incompressible Navier-Stokes. We justified the use of transfer learning for these specific test cases by studying the neural network mappings. We expect the applicability of transfer learning to generalise to other types of differential equations, provided that they possess similar smoothness properties, which is necessary for the error behaviour to be similar at different levels. In all test cases, fast convergence is obtained, leading to an accurate surrogate model already at a relatively low number of model runs. We compared our approach to MLSC. The transfer learning reduces the number of samples on the higher levels, and therefore significantly decreases the computational cost when medium/high fidelity samples are expensive to sample or when a surrogate is wanted with a high accuracy, which effectively increases the number of required levels. The conclusion is that our method outperforms MLSC when either the PDEs are sufficiently complex, or when a highly accurate surrogate model is required. For example, for the complex test case of a liquid sloshing in a tank, our approach leads to a computational cost savings of a factor of 2 compared to MLSC, when requiring medium to high accuracy. The resulting surrogate model can be directly used as a computationally inexpensive tool for uncertainty quantification.

Furthermore, the method can be applied to unsteady problems, but is not optimal for it in its current form. Either the neural networks have to be retrained for different time instances, or the temporal component should be added as an extra dimension to the convolutional layers and the output of the neural networks. Furthermore, it is not shown how the method scales for high-dimensional random spaces. Optimising the method proposed here for unsteady problems and high-dimensional random spaces is scheduled for future work.

\section*{Acknowledgements}
\noindent This work is part of the research programme ''SLING'' (Sloshing of Liquefied Natural Gas), which is (partly) financed by the Netherlands Organisation for Scientific Research (NWO).
\section*{References}
\bibliography{MyLibrary}
\end{document}